\documentclass[12pt]{amsart}
\usepackage{amscd}
\usepackage{amssymb,amsmath}
\usepackage[pagebackref]{hyperref}
\usepackage{mathrsfs}
\usepackage{graphicx}
\usepackage{enumitem}
\usepackage{romannum}
\usepackage{mathtools}
\usepackage[utf8]{inputenc}
\usepackage{listings}
\usepackage{cleveref}
\usepackage{amsmath}
\usepackage{tikz}
\usetikzlibrary{graphs}
\usetikzlibrary{decorations.markings}
\usepackage{lipsum,eso-pic,xcolor}
\usepackage{lineno}

\usepackage{tikz}
\usetikzlibrary{arrows,decorations.pathmorphing,backgrounds,positioning,fit}
\usetikzlibrary{positioning}
\hypersetup{
	colorlinks=true,
	linkcolor=blue,}

\usetikzlibrary{trees}
\usetikzlibrary{arrows}

\def\reg{\operatorname{reg}}
\def\deg{\operatorname{deg}}

\DeclarePairedDelimiter\ev{\langle}{\rangle}

\newcommand{\NN}{\mathbb{N}}

\newcommand{\N}{N}

\newcommand{\PP}{\mathcal{P}}

\newtheorem{lemma}{Lemma}[section]
\newtheorem{corollary}[lemma]{Corollary}
\newtheorem{theorem}[lemma]{Theorem}
\newtheorem{proposition}[lemma]{Proposition}
\newtheorem{definition}[lemma]{Definition}
\newtheorem{remark}[lemma]{Remark}
\newtheorem{example}[lemma]{Example}

\advance\headheight1.15pt

\begin{document}
	
	\pagenumbering{arabic}
	
	\title[Componentwise Linear Weighted Oriented Edge Ideals]{Componentwise linearity of powers of edge ideals of weighted oriented graphs} 
	
	\author[M. Kumar]{Manohar Kumar}
	\address{Department of Mathematics, Indian Institute of Technology
		Madras, Chennai, INDIA - 600036.}
	\email{manhar349@gmail.com}

    \author[J. Mondal]{Joydip Mondal$^*$}
	\address{Department of Mathematics, Indian Institute of Technology
		Kharagpur, West Bengal, INDIA - 721302.}
	\email{joydipmondal1999i@gmail.com}
	
	\author[R. Nanduri]{Ramakrishna Nanduri}
	\address{Department of Mathematics, Indian Institute of Technology
		Kharagpur, West Bengal, INDIA - 721302.}
  \email{nanduri@maths.iitkgp.ac.in}

        \thanks{$^*$ Supported CSIR-UGC PhD Fellowship, India}
	\thanks{AMS Classification 2020: 13D02, 05E40, 05E99, 13D45}

	%\keywords{Regularity of }
	\maketitle
\begin{abstract}		
In this paper, we study the componentwise linearity of powers of edge ideal of a weighted oriented graph $D$. We give a characterization for componentwise linearity of the edge ideal $I(D)$ in terms of forbidden subgraphs of $D$. If $D$ is house-free or complete $r$-partite, then the following statements are equivalent: 
\begin{center}
\begin{enumerate}
         \item $I(D)$ is componentwise linear;
         \item $I(D)$ is vertex splittable;
         \item $I(D)$ has linear quotient property;
         \item both $G$ and $H(I(D)_{\ev{2}})$ are co-chordal and $D_1,D_2,D_3,D_4$ as in Figure \ref{fig1}, are not induced subgraphs of $D$.
     \end{enumerate}
     \end{center}    
 Furthermore, if $D$ is a complete $r$-partite weighted oriented graph, then we show that 
 \begin{center}
 $I(D)^k \text{ is componentwise linear, for some } k\geq 2 \iff I(D) \text{ is componentwise linear}.$
 \end{center}
\end{abstract}

\section{Introduction}

Let $R=\mathbb{K}[x_1,\ldots,x_n]$ be the standard graded polynomial ring over a field $\mathbb{K}$ and $I\subset R$, a graded ideal with graded decomposition $I=\oplus_{d\geq 0} I_d$. For each $d$, let $I_{\ev{d}}:=(I_{d})$, the ideal generated by $I_d$. We say that $I$ is {\it componentwise linear} if $I_{\ev{d}}$ has a $d$-linear resolution for each $d\geq 0$, that is, $\reg(I_{\ev{d}})=d$. In 1999, Herzog and Hibi introduced these ideals and showed that the Stanley-Reisner ideal of a simplicial complex is componentwise linear if and only if its dual ideal is sequentially Cohen-Macaulay, see \cite{hh99}. Later many researchers studied various properties of the componentwise linear ideals, see \cite{hh99},  \cite{hi05}, \cite{nr15} \cite{hhm22}, \cite{hv22} etc. It is an active area of research in Commutative Algebra to study various algebraic and combinatorial properties of componentwise linear monomials ideals. R. Fr\"{o}berg in $1988$, in his celebrated theorem, known as Fr\"{o}berg's theorem, proved that the edge ideal of a finite simple graph is componentwise linear if and only if the graph is co-chordal. Recently in \cite{kns25}, the authors studied the componentwise linear edge ideals of vertex-weighted oriented graphs. A {\em (vertex) weighted oriented graph} is a triplet $D=(V(D), E(D), w)$, where $V(D)$ is the vertex set of $D$, $E(D)=\{(x,y) |\mbox{ there is a directed edge from vertex $x$ towards vertex $y$}\}$ is the {\it edge set} of $D$, and $w:V(D)\rightarrow \NN$ is a map, called weight function, i.e., a weight $w(x)$ is assigned to each vertex $x\in V(D)$. Note that $D$ has no multiple edges and no loop edges. Corresponding to a weighted oriented graph $D$, there is a simple graph $G$, called the underlying graph of $D$, such that $V(G):=V(D)$ and $E(G):=\{\{x,y\}| (x,y) \mbox { or } (y,x)\in E(D)\}$, i.e., $G$ is the simple graph without orientation and without weights in $D$. Let $D$ be a weighted oriented graph with $V(D)=\{x_1,\ldots,x_n\}$. Then the {\it edge ideal} of $D$, denoted by $I(D)$, is an ideal of $R$ defined as follows
\begin{equation*}
    I(D)=(x_ix_j^{w(x_j)}~|~(x_i,x_j)\in E(D)).
\end{equation*}

In \cite{kns25}, it was shown that if $I(D)$ is componentwise linear, then $\sqrt{I(D)}$ is componentwise linear. Also, for some classes of weighted oriented graphs $D$, the authors proved the following equivalence. 
\begin{center}
$I(D)$ is vertex splittable $\Leftrightarrow I(D)$ has linear quotients $\Leftrightarrow I(D)$ is componentwise linear. 
\end{center}
And they raised a question that whether the above equivalence hold for any $D$, \cite[Question 4.13]{kns25}. In this work, we prove the above equivalence for a significant classes of weighted oriented graphs and give a characterization for componentwise linearity. Also, we study the componentwise linearity of powers of edge ideals of weighted oriented graphs. 

The motivation for studying edge ideals of weighted oriented graphs comes from coding theory, in particular, in the study of Reed-Muller-type codes, see \cite{hlmrv19,prt19}. These ideals appear as the initial ideals of certain vanishing ideals in the theory of Reed-Muller-type codes. The study of weighted oriented edge ideals helps to obtain some properties of Reed-Muller codes easily. In particular, if $I(D)$ is Cohen-Macaulay, then vanishing ideal is Cohen-Macaulay. So far, a significant amount of research has been done concerning algebraic invariants and properties of these ideals (see \cite{gmsv18}, \cite{hlmrv19}, \cite{prt19}, \cite{x21}, \cite{kblo22}, \cite{bds23}, \cite{kn23}, \cite{kn25}). \par 

Let $D$ be a weighted oriented chordal graph or a $G^{\prime}$-free graph or a complete $r$-partite graph. Let $G$ be the underlying simple graph of $D$. Then we showed that the following statements are equivalent. 
\begin{enumerate}
         \item $I(D)$ is componentwise linear;
         \item $I(D)$ is vertex splittable;
         \item $I(D)$ has linear quotient property; 
         \item both $G$ and $H=H(I(D)_{\ev{2}})$ are co-chordal and $D_1,D_2,D_3,D_4$ are not induced subgraphs of $D$, 
     \end{enumerate}
where $G^{\prime}$ is the graph as in Figure \ref{fig8}, known as house graph and $D_i$'s are weighted oriented graphs as in Figure \ref{fig1} (Theorem \ref{theorem-3.4}, \ref{thm3}, \ref{compartite}). This class of house-free graphs includes chordal graphs, gap-free graphs, bipartite graphs, complete graphs, whisker graphs, $C_3$-free graphs, $C_4$-free graphs etc. Note that componentwise linearity depends on the characteristic of $\mathbb{K}$ (for example, see \cite[Remark 3]{reisner76}), but vertex splittable and linear quotient properties are characteristic independent.

The powers of componentwise linear ideals have been studying extensively for more than a decade. Many researchers focus on  when all powers of an ideal inherit componentwise linearity, their structural properties, see the recent survey in \cite{hv22}. The Herzog-Hibi-Ohsugi conjecture states that all powers of vertex cover ideals of chordal graphs are componentwise linear. This conjecture is still open. Powers of a componentwise linear ideal may fail to be componentwise linear. However, for the edge ideal $I(G)$ of simple graph $G$, we have $I(G)$ is componentwise linear if and only if  $I(G)^k$ is componentwise linear, for all $k\geq 2$ (\cite[Theorem 3.2]{hhz004}). It is not difficult to see that $I(D)^k$ need not be componentwise linear, if $I(D)$ is componentwise linear, for a weighted oriented graph $D$ (see Example \ref{power}). If $D$ is a complete $r$-partite $K_{n_1,\dots,n_r}$, then we show that 
$$I(D)^k \text{ is componentwise linear, for some } k\geq 2 \iff I(D) \text{ is componentwise linear}$$
(see Theorem \ref{thm4}). Furthermore, we show that 
$$I(D)^2 \text{ has linear quotients} \iff I(D) \text{ has linear quotients.}$$ 
 As a consequence, if $I(D)$ is not componentwise linear, then $I(D)^k$ is not componentwise linear, for all $k\geq 2$. The key lemmas to prove our main results is that if $D$ is complete $r$-partite, then $|V^+|\leq (r-1)$, where $V^+=\{x\in V(D)~|~ w(x) \geq 2 \text{ and } x \text{ is not a source}\}$ (Lemmas \ref{lemma-4.6}, \ref{verpartite}). 
 
Now, we describe the organization of the paper. In Section \ref{secpreli}, we recall definitions, notions, and results associated with our work. In section \ref{secvertexsplit}, we have shown a characterization for componentwise linearity of edge ideals of chordal graphs, house-free graphs, and complete $r$-partite graphs. In Section \ref{secpowervertexsplit}, we investigate componentwise linearity of powers of house-free weighted oriented graphs. Finally in Section \ref{sec5}, we study the componentwise linearity of powers of complete $r$-partite weighted oriented graphs.

\section{Preliminaries}\label{secpreli}

In this section, we recall some necessary prerequisites, which are used to describe our work and establish our results.\par 

Let $D=(V(D), E(D), w)$ be a (vertex) weighted oriented graph with $V(D)=\{x_1,\ldots,x_n\}$ and underlying simple graph $G$. An edge $e\in E(D)$ is an ordered pair $e=(x_i,x_j)$ with $x_i,x_j\in V(D)$ and the orientation of $e$ is from the vertex $x_i$ to the vertex $x_j$. Also, $w:V(D)\rightarrow \NN$ is the weight function which assigns a weight for each vertex of $D$. Let $\mathfrak{m}=(x_1,\ldots,x_n)$ be the homogeneous maximal ideal of $R=\mathbb{K}[x_1,\ldots,x_n]$, where $\mathbb{K}$ is a field. We write $V^+(D):= \{x\in V(D): w(x)>1 ~~\text{and}~~x ~\text{is not a  source vertex of $D$}\}$, in short $V^{+}(D)$ is denoted by $V^+$. For a vertex $x\in V(D)$, its {\it outer neighbourhood} is defined as 
  $N_{D}^{+}(x):= \{y \in V(D) | (x,y)\in E(D)\}$, its {\it inner neighbourhood} is defined as $N_{D}^{-}(x):= \{z \in V(D) | (z,x)\in E(D)\}$, and $N_D(x):=N_{D}^{+}(x) \cup N_{D}^{-}(x)=N_{G}(x)$. Also, we denote $N_D[x]:=N_D(x) \cup \{x\}=N_{G}[x]$. A vertex $x\in V(D)$ is called a {\it source} if $N_{D}^{-}(x)=\emptyset$ and $x$ is called a {\it sink} if $ N_{D}^{+}(x) = \emptyset $. If $x\in V(D)$ is a source, then we may assume $w(x)=1$ as it would not affect the ideal $I(D)$. The \textit{degree} of a vertex $x$ in $G$ is defined as $\text{deg}_G(x) :=|N_G(x)|$. For a subset $A\subseteq V(D)$, we denote the induced subgraph of $D$ on $A$ by $D[A]$. If $x\in V(D)$ is a vertex, then we write $D\setminus x$ to denote the graph $D[V(D)\setminus \{x\}]$.
%Also, we denote $H(I(D)_{\ev 2})$ as a subgraph (need not be induced subgraph) $H$ of the simple graph $G$ such that $I(H)=I(D)_{\ev 2}$, where $G$ is the underlying simple graph of $D$. \par
  
 Now, let us discuss some notations and definitions regarding simple graphs. A graph $G$ is said to be a {\it complete graph} if there is an edge between each pair of vertices of $G$ and a complete graph on $n$ vertices is denoted by $K_{n}$. A {\it clique} of a graph $G$ is a set $A$ of vertices of $G$ such that $G[A]$ is a complete graph. A vertex $x\in V(G)$ is said to be a {\it simplicial vertex} of $G$ if $G[N_{G}[x]]$ is a complete graph. A subset $B \subseteq V(G)$ is called an {\it independent set} if no two vertices in $B$ are adjacent in $G$. A set of vertices $C$ of $G$ is called a {\it vertex cover} of $G$ if $C\cap e\neq \emptyset$ for all $e\in E(G)$. A {\it minimal vertex cover} is a vertex cover that is minimal with respect to inclusion. The \textit{complement} of a simple graph $G$, denoted by $G^c$, is the simple graph such that $V(G^c)=V(G)$ and $E(G^c)=\{\{u,v\}\mid \{u,v\}\not\in E(G)\}$. A cycle of length $n$, denoted by $C_n$, is a connected graph such that every vertex of $C_n$ has degree two. A simple graph $G$ is called {\it chordal} if $G$ has no induced cycle of length greater than three. A graph $G$ is called {\it co-chordal} if $G^c$ is chordal. A graph $G$ is said to be {\it bipartite} if it has no induced odd cycle. A weighted oriented graph $D$ is said to be a cycle or complete or chordal or co-chordal or bipartite if its underlying simple graph is so. We say that a finite simple graph $G$ is house-free graph, if $G$ does not contain $G^{\prime}$ as an induced subgraph of $G$, where $G^{\prime}$ is as shown in the Figure \ref{fig8}. One can define a weighted oriented graph $D$ is house-free if its underlying simple graph is house-free. 

% \begin{definition}{\rm For any homogeneous ideal $I$ in $R$, the Castelnuovo-Mumford regularity (or simply regularity), denoted by $\reg(I)$, is defined as follows 
%  \begin{align*}
%  \reg(I) &=  \max \{j - i \mid \beta_{i,j}(I) \neq 0\} \\
%          &= \max\{j+i \mid H_{\m}^i(I)_j \neq 0\},   
%  \end{align*}
%   where $\beta_{i,j}(I)$ is the $(i,j)^{th}$ graded Betti number of $I$ and $H_{\m}^i(I)_j$ denotes the $j^{th}$ graded component of the $i^{t h}$ local cohomology module $H_{\m}^i(I)$. 
%   }
% \end{definition}
 
\begin{definition}{\rm
   A homogeneous ideal $I$ of $R$ is componentwise linear if for each $d\geq 0$, the ideal generated by $d^{th}$ homogeneous component $I_d$, denoted by $I_{\ev{d}}$, has $d$-linear resolution over $R$, i.e. $\reg(I_{\langle d \rangle})=d$, where $\reg(-)$ denote the Castelnuovo-Mumford regularity. Note that if $I$ is equigenerated, then linearity and componentwise linearity are equivalent.
   }
\end{definition}
 
\begin{definition}{\rm
    A monomial ideal $I \subseteq R$ has linear quotient property if there exists an order $ u_1 < \cdots < u_m$ on the minimal monomial generating set $\mathcal{G}(I)=\{u_1,\ldots,u_m\}$ of $I$ such that the colon ideal $((u_1, \ldots, u_{i-1}) : u_i)$ is generated by a subset of the variables, for $i = 2, \ldots, m$.
    }
\end{definition}

Now, let us state the well-known Fr\"{o}berg's theorem \cite{f90}, which has been used frequently in this paper.\medskip

\noindent{\bf Fr\"{o}berg's Theorem.} Let $G$ be a simple graph and $I(G)$ be its edge ideal. Then $I(G)$ has linear resolution if and only if $G^c$ is chordal.

\begin{definition}\label{spl1}{\rm 
A monomial ideal $I\subseteq R=\mathbb{K}[X]$ is called {\it vertex splittable} if it can be obtained by the following recursive procedure.
\begin{enumerate}
    \item If $v$ is a monomial and $I=(v)$, $I=(0)$ or $I=R$, then $I$ is vertex splittable.
    \item If there is a variable $x$ in $R$ and vertex splittable ideals $I_1$ and $I_2$ in $\mathbb{K}[X\setminus \{x\}]$ such that $I=xI_1+I_2$, $I_2 \subseteq I_1 $ and $\mathcal{G}(I)= \mathcal{G}(xI_1)\sqcup \mathcal{G}(I_2)$, then $I$ is a vertex splittable. For $I=xI_1+I_2$, the variable $x$ is said to be a {\it splitting variable} for $I$.
\end{enumerate}
}
\end{definition}

%\noindent 
%Now recall the following lemma and its proof follows from \cite[Lemma 7]{kmns13}. 
%\begin{lemma}\label{simplicial}
%Every chordal graph $G$ has a simplicial vertex. Moreover, if $G$ is a non-complete chordal graph, then it has two non-adjacent simplicial vertices.
%\end{lemma}

\begin{theorem}\cite[Corollary 3.3]{kns25}\label{corollary-2.7}
    Let $D$ be a weighted oriented graph. Suppose $I(D)$ is componentwise linear. Then $D_1, D_2, D_3$ and $D_4$ as in Figure \ref{fig1} cannot be induced subgraphs of $D$.
\end{theorem}

\begin{lemma}\cite[lemma 4.1(1)]{kns25}\label{lemma-2.8}
    Let $I \subseteq R$ be a monomial ideal. If $I$ is vertex splittable, then $(I,x_{i_1},\ldots,x_{i_r})$ is vertex splittable, where none of the variables $x_{i_{j}}$ are in \text{supp}$(\mathcal{G}(I))$. 
\end{lemma}

\begin{lemma}\cite[Lemma 2.1]{jz10}\label{increasing-admissible}
Let $I \subseteq R$ be a monomial ideal with linear quotients. Then there is a degree increasing admissible
order of $\mathcal{G}(I)$ such that $I$ has linear quotients with respect to this degree increasing order.
\end{lemma}

%\begin{theorem}\cite[corollary 3.12]{bed24}\label{powerlinear}
%  Let $G$ be a simple graph and $I(G)$ denote its edge ideal. Then $I(G)$ has linear quotients if and only if $I(G)^k$ has linear quotients for all $k \geq 1$.  
%\end{theorem}

\begin{theorem} \cite[Theorem 3.7]{kns25}\label{thmcochordal}
	Let $D=(V(D), E(D), w)$ be a weighted oriented graph with the underlying simple graph $G$. If $I(D)$ is componentwise linear, then $G^c$ is chordal.
\end{theorem}

\begin{lemma}\label{inducedpower}
    Let $D$ be a weighted oriented graph. If $I(D)^k$ is componentwise linear, then $I(H)^k$ is also componentwise linear for any induced subgraph $H$ of $D$.
\end{lemma}
\begin{proof}
    This follows from \cite[Proposition-3.1(2)]{kns25}.
\end{proof}
% \begin{figure}[!h]
%     \centering
%      \includegraphics[width=0.4 \textwidth]{figure5.png} 
%        \caption{A weighted oriented cycle $\overline{D}$ with $I(\overline{D})$ componentwise linear}
%     \label{fig5}
% \end{figure}

% \begin{theorem}\cite[Theorem 4.9]{kns25}\label{thmdbar}
%      Let $D=(V(D), E(D), w)$ be a weighted oriented chordal graph, and the graph $\overline{D}$ in Figure \ref{fig5} is not an induced subgraph of $D$. Then, the following statements are equivalent.
%      \begin{enumerate}
%          \item $I(D)$ is componentwise linear;
%          \item $I(D)$ is vertex splittable;
%          \item $I(D)$ has linear quotient property.
%      \end{enumerate}
% \end{theorem}

\section{Componentwise linear $G'$-free weighted oriented graphs} \label{secvertexsplit}

\begin{figure}[!h]
    \centering
     \includegraphics[width=0.6 \textwidth]{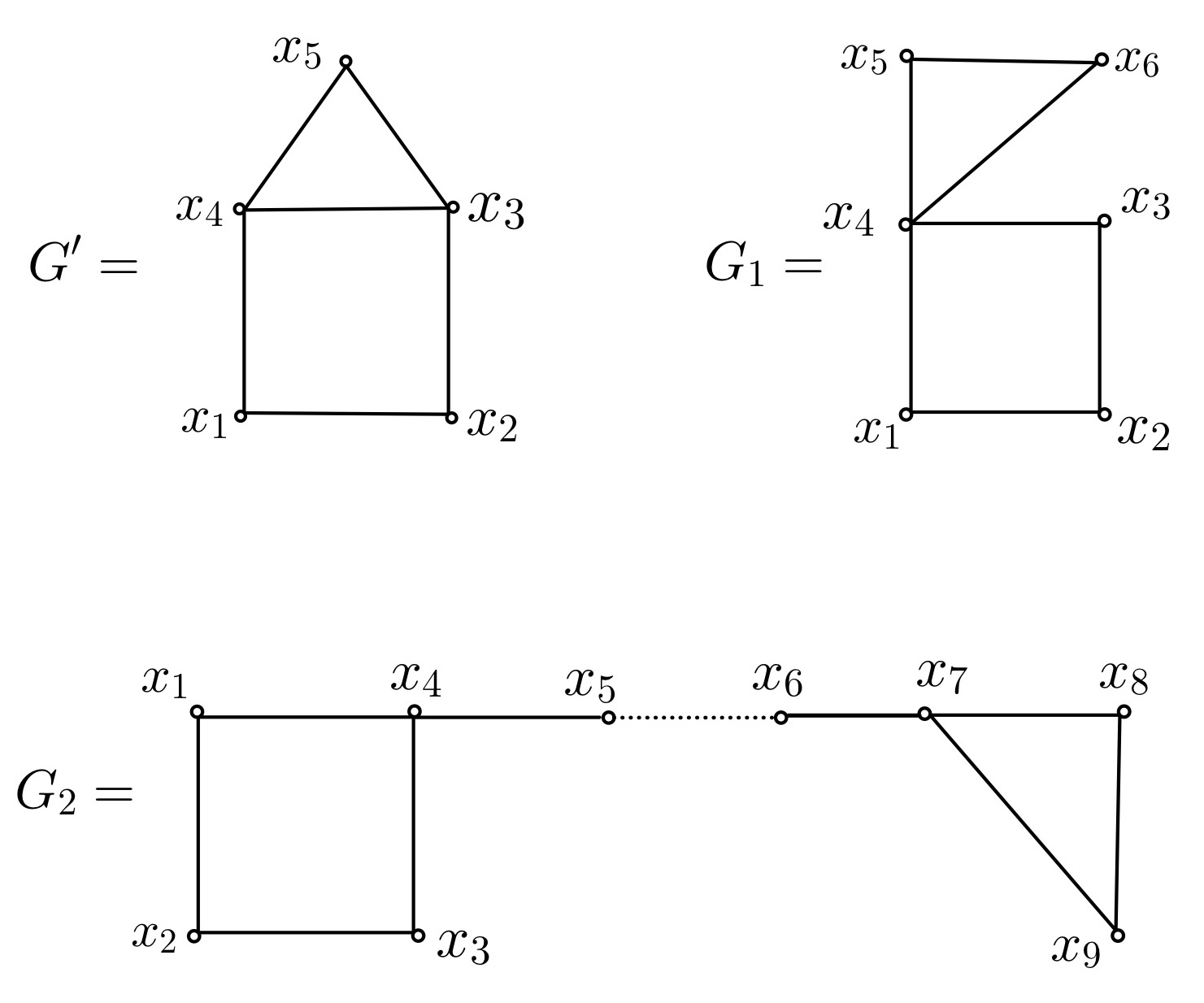} 
    \caption{$G^{\prime}$ is the house graph, $G_1$ is a $3$-cycle and a $4$-cycle sharing a vertex and $G_2$ is a $3$-cycle and a $4$-cycle connected by a path}
    \label{fig8}
\end{figure}

In this section, we give a characterization for componentwise linearity of $I(D)$, where $D$ is house-free (that is, $G^{\prime}$-free) or complete $r$-partite, where $G^{\prime}$ is as in Figure \ref{fig8}. 

\begin{lemma}\label{lemma-3.1}
    Let $D$ be a weighted oriented graph. If $D_i$ are not induced subgraphs of $D$ for $i\in \{1,\ldots,4\}$,  then  $|\N_D^+(x)\cap V^+|\leq 1$, for all $x\in V(D)$.
\end{lemma}
\begin{proof}
   Let $D$ be a weighted oriented subgraph and $D_i$ are not induced subgraphs of $D$ for all $i\in\{1,\ldots,4\}$. Suppose there exists a vertex $x\in V(D)$ such that $|N_D^+(x)\cap V^+|> 1$. Let $x_1,x_2\in N_D^+(x)\cap V^+$. If we consider the induced weighted subgraph $D'=D[\{x,x_1,x_2\}]$, then we have three possible edge set of $D$
  \begin{enumerate}
      \item $ E(D')=\{(x,x_1),(x,x_2)\}$. This implies $D'=D_2$.
      \item $E(D')=\{(x,x_1),(x,x_2),(x_2,x_1)\}$. This implies $D'=D_4$.
      \item $E(D')=\{(x,x_1),(x,x_2),(x_1,x_2)\}$. This implies $D'=D_4$
  \end{enumerate}
  In all the cases, we get a contradiction. Thus $|N_D^+(x)\cap V^+|\leq 1$, for all $x\in V(D)$.
\end{proof}

\begin{figure}[!h]
    \centering
     \includegraphics[width=0.4 \textwidth]{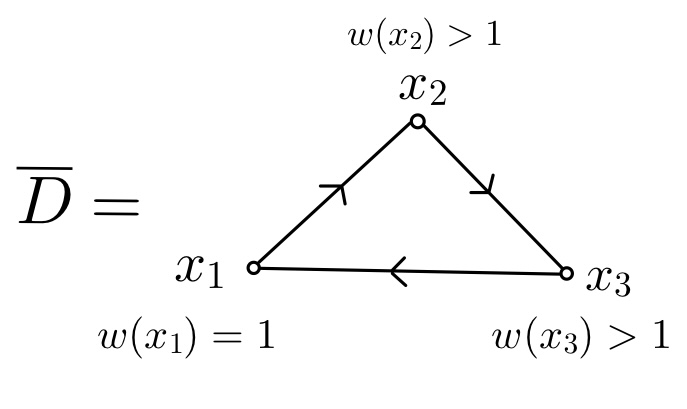} 
       \caption{A weighted oriented cycle $\overline{D}$ with $I(\overline{D})$ componentwise linear}
    \label{fig5}
\end{figure}

\begin{figure}[!h]
    \centering
     \includegraphics[width=1.0 \textwidth]{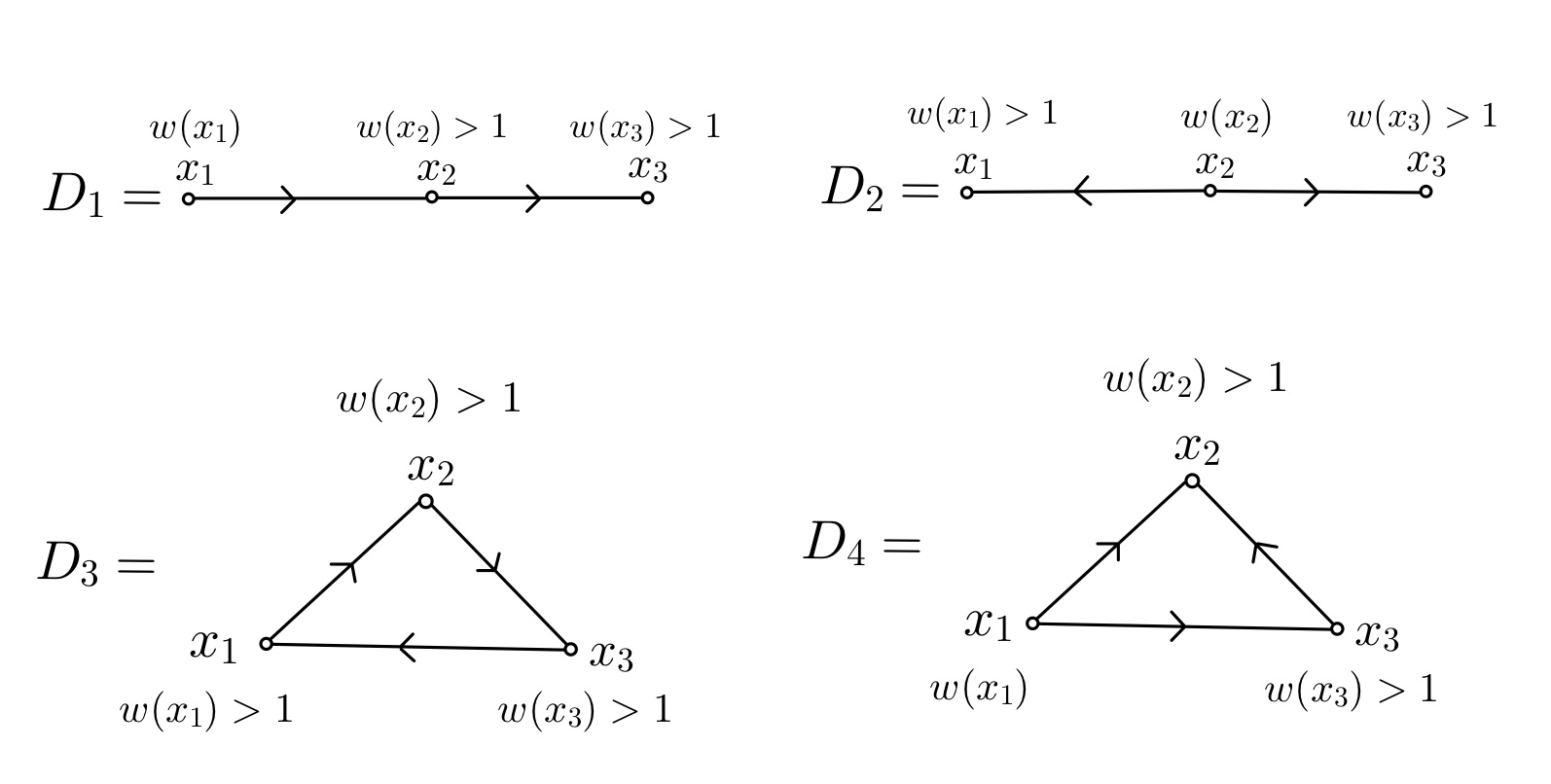} 
    \caption{$D_1, D_2, D_3, D_4$ can not be induced subgraphs of a weighted oriented graph $D$ whose edge ideal is componentwise linear.}
    \label{fig1}
\end{figure}

Below we prove a characterization for componentwise linearity of edge ideals of chordal graphs. 

\begin{theorem}\label{theorem-3.4}
     Let $D=(V(D), E(D), w)$ be a weighted oriented chordal graph and $G$ denote its underlying simple graph. Then, the following statements are equivalent.
     \begin{enumerate}
         \item $I(D)$ is componentwise linear;
         \item $I(D)$ is vertex splittable;
         \item $I(D)$ has linear quotient property; 
         \item both $G$ and $H=H(I(D)_{\ev{2}})$ are co-chordal and $D_i$ are not induced subgraphs of $D$ for $i=1,\ldots,4$.
     \end{enumerate}
\end{theorem}      
\begin{proof}
$(2)\implies (3) \implies (1)$: Clear. 

\begin{figure}[!h]
    \centering
     \includegraphics[width=0.5 \textwidth]{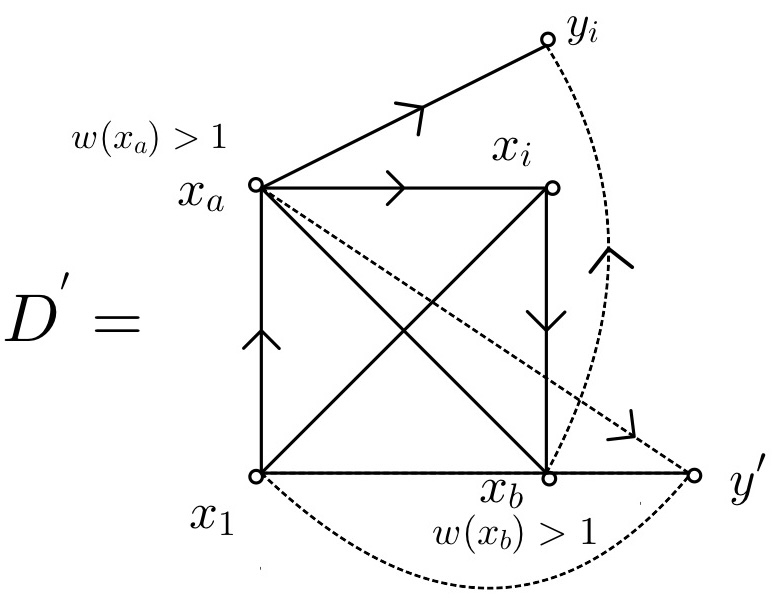} 
    \caption{A non co-chordal induced subgraph of $D$}
    \label{fig7}
\end{figure}

\noindent $(1)\implies (4)$: Let $I(D)$ is componentwise linear then by Theorem \ref{corollary-2.7} we get $D_i$ are not induced subgraphs of $D$ for $i\in \{1,\ldots,4\}$. Since $I(D)$ is componentwise linear then $I(H)=I(D)_{\ev{2}}$ has linear resolution. Therefore, by { Fr\"{o}berg's Theorem} we get $H$ is co-chordal. By Theorem \ref{thmcochordal}, $G$ is co-chordal. \\
   \noindent $(4)\implies (2)$: Let $n$ be the number of vertices of $D$. We may assume that $n\geq 2$ and $D$ is not an empty graph. Note that both $G$ and $G^c$ are chordal. Thus, $G$ is a split graph. Consequently, we can partition $V(G)$ as $V(G)=A\cup B$, where $G[A]$ is a maximal clique and $B$ is an independent set of $G$. Since $G$ and $H$ both are co-chordal, $G\setminus z$ and $H\setminus z$ are also co-chordal for $z\in V(D)$. Also, it is easy to verify that $H\setminus z=H(I(D\setminus z)_{\ev{2}})$ for $z\in V(D)$. Let $A=\{x_1,\ldots, x_r\}$ and $B=\{y_1,\ldots,y_s\}$. Note that $N_{G}(x_i)$ is a minimal vertex cover of $G$ for all $1\leq i\leq r$. Suppose $\vert A\cap V^{+}\vert> 2$, then there exist $x,w,z \in A\cap V^+$ such that the induced subgraph $D[\{x,w,z\}]$ is either $D_3$ or $D_4$ which is a contradiction. Therefore $\vert A\cap V^{+}\vert\leq 2$. We proceed by induction on $n$. The case $n=2$ is obvious. To complete the proof, it is sufficient to prove the following cases: \\
\textnormal{(i)} \((4) \Rightarrow (2)\) when \(r \leq 3\); \\
\textnormal{(ii)} \((4) \Rightarrow (2)\) when \(r \geq 4\). \\
  \textbf{\underline{Proof of (i):}} Assume $r \leq 3$. If $r=2$ or ( $r= 3$ and $\vert A\cap V^{+}\vert\leq 1$), then it is easy to see that $\overline{D}$ is not an induced subgraph of $D$. Then by \cite[Theorem 4.9]{kns25}, the required implication holds. Now, assume that $r=3$ and $|A\cap V^+|=2$. Without loss of generality, assume that $x_2,x_3 \in V^+$ and $w(x_1)=1$. Assume $(x_2,x_1)$ and $(x_3,x_1)\in E(D)$ then, we can write
   $I(D)=x_1I_1+I(D\setminus x_1)$,
   where $$I_1=\left(z^{w(y)}\mid z\in N_{D}^{+}(x_1)\cap V^{+}\right )+\left (z\mid z\in (N_{D}^{+}(x_1)\setminus V^{+}) \cup N_{D}^{-}(x_1)\right).$$ Thus by Lemma \ref{lemma-2.8} and \ref{lemma-3.1}, $I_1$ is vertex splittable. Since $(A\setminus \{x_1\})\subseteq I_1$, we get that $I(D\setminus x_1)\subseteq I_1$. Now by the induction hypothesis, $I(D\setminus x_1)$ is vertex splittable. Thus $I(D)$ is vertex splittable.\\
   Assume $(x_1,x_2)\in E(D)$, then $x_2\in N_D^+(x_1)\cap V^+$. Then using Lemma \ref{lemma-3.1} we have $(x_3,x_1)\in E(D)$. Now, we will take two possible cases:\\
   \textbf{Case-1:}
    Let $(x_2,x_3)\in E(D)$. We can write $I(D)=x_1I_1+I(D\setminus x_1)$, where $$I_1=\left(z^{w(z)}\mid z\in N_{D}^{+}(x_1)\cap V^{+}\right)+\left(z\mid z\in (N_{D}^{+}(x_1)\setminus V^{+})\cup N_{D}^{-}(x_1)\right).$$ Using Lemma \ref{lemma-3.1} we have $|N_D^+(x_1)~\cap~ V^+|\leq 1$. Since $x_2\in N_D^+(x_1)\cap V^+$ this implies that $I_1=(x_2^{w{(x_2)}})+ \left(z\mid z\in (N_{D}^{+}(x_1)\setminus V^{+}) \cup N_{D}^{-}(x_1)\right ).$
    By induction, $I(D\setminus x_1)$ is vertex splittable, and by Lemma \ref{lemma-2.8}, $I_1$ is vertex splittable. If $N_D^+(x_2)\setminus N_D(x_1)=\emptyset$, then $I(D\setminus x_1)\subseteq I_1 $ which yields that $I(D)$ is vertex splittable. \\
    Assume $N_D^+(x_2)\setminus N_D(x_1)\neq \emptyset$. Let $y\in (N_D^+(x_2)\setminus N_D(x_1))\cap B$. Then $w(y)=1$ as $y\notin N_D(x_1)$ and $D_1$ can not be an induced subgraph of $D$. If possible let $y_j\in N_D(x_1)\cap B$ for some $j$. Suppose $w(y_j)>1$ then $\{y_j,x_1\}$ and $\{x_2,y\}$  form an induced matching in $H$. This is a contradiction as $H$ is co-chordal. This implies that $w(y_j)=1$. Since $\{x_3,x_1\}$ and $\{x_2,y\}$ can not form an induced matching in $H$, we should have $\{x_3,y\}\in E(H)$. Also $\{y_j,x_1\}, \{x_2,y\}$ cannot form an induced matching in $H$ then we have $\{x_2,y_j\}\in E(H)$. We have $G[\{x_1,x_2,x_3\}]$ is a maximal clique of $G$ then $(x_3,y_j)\notin E(D)$. This implies that $G[\{x_1,y_j,x_2,y,x_3\}]$ will form an induced $5$-cycle in $H$. This is a contradiction because $H$ is co-chordal. Therefore, $y_j\notin N_D(x_1)\cap B$ for all $j$. Thus we can write $I(D)$ as follows, $I(D)=y(x_2,x_3)+I(D\setminus y)$. Since $I(D\setminus y)\subseteq (x_2,x_3)$ and by induction hypothesis $I(D\setminus y)$ is vertex splittable. Hence $I(D)$ is vertex splittable. \\
    \textbf{Case-2:}
    Let$( x_3,x_2)\in E(D)$. Since $x_3\in V^+$ and  $(x_3,x_1) \in E(D)$ then there exists a $y'\in B$ such that $(y',x_3)\in E(D)$ and $w(y')=1$. Since $D_1$ and $D_4$ are not induced subgraphs of $D$ then we have  $(x_2,y')\in E(D)$. Now in $H$, either $\{x_2,y'\}$ and $\{x_3,x_1\}$ will form an induced matching or $(x_1,y')\in E(D)$. If  $\{x_2,y'\}$ and $\{x_3,x_1\}$ form an induced matching in $H$, which contradict that $H$ is co-chordal. If  $(x_1,y')\in E(D)$, then we get  $G[\{x_1,x_2,x_3,y'\}]$ is a maximal clique of $G$ which is a contradiction because $G[\{x_1,x_2,x_3\}]$ is a maximal clique of $G$. In both cases, we arrive a contradiction and hence $(x_3,x_2)\notin E(D)$ i.e., $(x_2,x_3)\in E(D)$. Therefore Case-2 never arise. \\ 
\textbf{\underline{Proof of (ii):}} 
    Let us assume that $r\geq 4$. Since $D_3$, $D_4$  cannot be induced subgraphs of $D$, we have $\vert A\cap V^{+}\vert\leq 2$. If $\vert A~\cap V^{+}\vert\leq 1$ then $\overline{D}$ is not an induced subgraph of $D$, then by Theorem \cite[Theorem 4.9]{kns25}, the required implication holds. Now, assume  $|A\cap V^+|=2$. Let $x_a, x_b \in V^+$. Now, we choose a vertex $x_1\in A~\cap~(V\setminus V^{+})$, with the highest degree, i.e., $\deg_{G}(x_1)\geq \deg_{G}(x_i)$ for all $x_i\in A\cap(V\setminus V^{+})$.  If $(x_a,x_1),(x_b,x_1)\in E(D)$ then, we can write
   $I(D)=x_1I_1+I(D\setminus x_1)$,
   where $I_1=\left(z^{w(z)}\mid z\in N_{D}^{+}(x_1)\cap V^{+}\right )+\left (z\mid z\in (N_{D}^{+}(x_1)\setminus V^{+}) \cup N_{D}^{-}(x_1)\right).$ 
   Using Lemma \ref{lemma-3.1} we get $\vert N_{D}^{+}(x_1)\cap V^{+}\vert\leq 1$, then by Lemma \ref{lemma-2.8}, $I_1$ is vertex splittable. We have $I(D\setminus x_1)\subseteq I_1$ as $(A\setminus \{x_1\})\subseteq I_1$. By the induction hypothesis, $I(D\setminus x_1)$ is vertex splittable, and hence $I(D)$ is vertex splittable.\\
   Assume $(x_1,x_a)\in E(D)$ then $x_a\in N_D^+(x_1)\cap V^+$. Then using Lemma \ref{lemma-3.1} we have $(x_b,x_1)\in E(D)$.
  We can write $I(D)$ as follows,
    $I(D)=x_1I_1+I(D\setminus x_1)$, where $I_1=(x_a^{w(x_a)})+ (z\mid z\in (N_D^+(x_1)\setminus x_a)\cup N_D^-(x_1))$.
    By induction, $I(D\setminus x_1)$ is vertex splittable, and $I_1$ is vertex splittable by Lemma \ref{lemma-2.8} and \ref{lemma-3.1}. If $N_D^+(x_a)\setminus N_D(x_1)=\emptyset$, then $N_{D}^{+}(x_a)\subseteq N_{D}(x_1)$. This implies that $I(D\setminus x_1)\subseteq I_1$ which yields that $I(D)$ is vertex splittable, as required. Suppose $N_D^+(x_a)\setminus N_D(x_1)\neq \emptyset$ then there exist a $y\in (N_D^+(x_a)\setminus N_D(x_1))\cap\ B$. Since $D_1$ is not an induced subgraph of $D$, then $w(y)=1$. Let $U=\{1,\ldots ,r\}\setminus\{1,a,b\}$. Now, we will consider two possible cases:\\
    \textbf{Case-1:}
    Assume  $(x_i,x_a)\in E(D)$ for some $i\in U$. Since $\{x_1,x_i\}$ and $\{x_a,y\}$ cannot form an induced matching in $H$, we should have  $\{y,x_i\}\in E(G)$ as $\{y,x_1\}\not\in E(G)$. Since $\deg_{G}(x_1)\geq \deg_{G}(x_i)$, there exists $y''\in N_{D}(x_1)\cap B$ such that $\{y'',x_i\}\not\in E(G)$. Again, since $D_2$ and $D_4$ cannot be induced subgraphs of $D$, we have $x_1y''\in I(D)$. Now in $H$, $\{x_1,y''\}$ and $\{x_a,y\}$ cannot form an induced matching, which gives that $w(y'')=1$ and $\{x_a,y''\}\in E(H)$. Thus the induced subgraph $G[\{x_1,y'',x_a,y,x_i\}]$ is an induced $5$-cycle in $H$. But this is a contradiction because $H$ is co-chordal. Therefore $N_D^+(x_a)\setminus N_D(x_1) = \emptyset$.  \\
     \textbf{Case-2:}
   Assume $(x_a,x_i)\in E(D)$ for all  $i\in U$. Suppose $(x_b,x_i)\in E(D)$ for some  $i\in U$. Then we can write as follows $I(D)=x_iI_1+I(D\setminus x_i)$,
   where $$I_1=\left(z^{w(z)}\mid z\in N_{D}^{+}(x_i)\cap V^{+}\right)+ (z\mid z\in  \left (N_{D}^{+}(x_i)\setminus V^{+}) \cup N_{D}^{-}(x_i)\right).$$ Then by Lemma \ref{lemma-2.8} and Lemma \ref{lemma-3.1}, $I_1$ is vertex splittable. Since $(A\setminus \{x_i\}) \subseteq I_1$, we have $I(D\setminus x_i)\subseteq I_1$. By the induction hypothesis, $I(D\setminus x_i)$ is vertex splittable, and hence $I(D)$ is vertex splittable.\\
    Suppose $(x_i,x_b)\in E(D)$ for all  $i\in U$. If for some $i\in U$, $N_D^+(x_b)\setminus N_D(x_i)= \emptyset$ i.e., $N_D^+(x_b)\subseteq N_D(x_i)$ then we can write $I(D)$ as follows,
    $I(D)=x_iI_1+I(D\setminus x_i)$, where $$I_1=\left(z^{w(z)}\mid z\in N_{D}^{+}(x_i)\cap V^{+}\right)+\left(z\mid z\in (N_{D}^{+}(x_i)\setminus V^{+}) \cup N_{D}^{-}(x_i)\right)~~ \text{and}~~ I_2=I(D\setminus x_i).$$ Using Lemma \ref{lemma-3.1}, we have $|N_D^+(x_i)~\cap~ V^+|\leq 1$. Since $x_b\in N_D^+(x_i)~\cap ~V^+$ this implies that $I_1=\left(x_b^{w(x_b)}\right)+ \left(z\mid z\in (N_{D}^{+}(x_i)\setminus V^{+}) \cup N_{D}^{-}(x_i)\right )$.
    By induction, $I_2$ is vertex splittable, and also $I_1$ is vertex splittable by Lemma \ref{lemma-2.8}. Since $N_D^+(x_b)\setminus N_D(x_i)=\emptyset$ then $I(D\setminus x_i)\subseteq I_1 $. Then $I(D)$ is vertex splittable.\\
 Assume $N_D^+(x_b)\setminus N_D(x_i)\neq \emptyset$, for all $i\in U$. Then there is an element  $y'\in (N_D^+(x_b)\setminus N_D(x_i)) \cap\ B$, for all $i$. Since $D_1$ is not an induced subgraph of $D$, then $w(y')=1$. If $r > 4$, then for two such vertices $j_1,j_2\in U$ we get an induced matching $\{x_{j_1},x_{j_2}\}$ and $\{x_b,y'\}$ in $H$ which is a contradiction because $H$ is co-chordal. Therefore $r\leq4$. Let $A=\{x_1,x_a,x_i,x_b\}$ and $y'\in (N_D^+(x_b)\setminus N_D(x_i)) \cap\ B$. Since $\{x_a,x_i\}$ and $\{x_b,y'\}$ cannot form an induced matching in $H$, we should have  $\{x_a,y'\}\in E(H)$ because $\{x_i,x_b\},\{y',x_i\},\{x_a,x_b\}\notin E(H)$. Thus $\{y',x_1\}\in E(H)$ otherwise $G[\{x_a,x_i,x_1,x_b,y'\}]$ is an induced $5$-cycle in $H$ which is absurd. Also $\{x_b,x_1\}$ and $\{x_a,y\}$ cannot form an induced matching in $H$ then we have $\{x_b,y\}\in E(H)$ as $\{x_1,y\}\notin E(H)$. Let $D'=D[\{y,x_a,x_i,x_1,x_b,y'\}]$ as in Figure \ref{fig7} be an induced subgraph of $D$, and let $H'=H[\{y,x_a,x_i,x_1,x_b,y'\}]$. Then we have $H^{\prime}$ is a subgraph of the underlying simple graph of $D^{\prime}$. One can easily check that $H'$ is not co-chordal. Since  $H'$ is an induced subgraph of $H$, this gives that $H$ is not co-chordal. This is a contradiction. Therefore $N_D^+(x_b)\setminus N_D(x_i) = \emptyset$ for some $i$. This implies that $N_D^+(x_b)\subseteq N_D(x_i) $. Then $I(D\setminus x_i)\subseteq I_1$. Therefore we can write $I(D)$ as follows,
    $I(D)=x_iI_1+I(D\setminus x_i)$, where $$I_1=\left(z^{w(z)}\mid z\in N_{D}^{+}(x_i)\cap V^{+}\right)+\left(z\mid z\in (N_{D}^{+}(x_i)\setminus V^{+}) \cup N_{D}^{-}(x_i)\right).$$
   Using Lemma \ref{lemma-3.1} we have $|N_D^+(x_i)~\cap~ V^+|\leq 1$. Since $x_b\in N_D^+(x_i)~\cap ~V^+$ this implies that $I_1=(x_b^{w(x_b)})+\left(z\mid z\in (N_{D}^{+}(x_i)\setminus V^{+}) \cup N_{D}^{-}(x_i)\right )$.
    By induction, $I(D\setminus x_i)$ is vertex splittable, and $I_1$ is vertex splittable by Lemma \ref{lemma-2.8}. Hence $I(D)$ is vertex splittable, as required.
\end{proof}

Now, we give a characterization for componentwise linearity of edge ideals of house-free (i.e., $G^{\prime}$-free) graphs. This class of graphs includes chordal graphs, gap-free graphs, bipartite graphs, complete graphs, whisker graphs, $C_3$-free graphs, $C_4$-free graphs etc. 

\begin{theorem} \label{thm3}
 Let $D$ be a weighted oriented graph and $G$ its underlying simple graph such that $G$ is $G^{\prime}$-free (that is, house-free),  where $G'$ is as in Figure \ref{fig8}. Then the following statements are equivalent.
  \begin{enumerate}
         \item $I(D)$ is componentwise linear;
         \item $I(D)$ is vertex splittable;
         \item $I(D)$ has linear quotient property; 
         \item both $G$ and $H=H(I(D)_{\ev{2}})$ are co-chordal and $D_1,D_2,D_3,D_4$ are not induced subgraphs of $D$.
     \end{enumerate}
\end{theorem}
\begin{proof}
 $(2) \implies (3) \implies (1)$ is easy. \\ 
\noindent $(1)\implies (4)$: Let $I(D)$ be componentwise linear. Then by Theorem \ref{corollary-2.7} we get $D_i$ are not induced subgraph of $D$ for $i\in \{1,\ldots,4\}$. Also by Theorem  \ref{thmcochordal}, we have $G$ is co-chordal. Since $I(D)$ is componentwise linear, we have  $H=H(I(D)_{\ev{2}})$ has linear resolution. Therefore, by { Fr\"{o}berg's Theorem} we get $H$ is co-chordal. This proves (4). \\
\noindent $(4)\implies (2)$: Let $G$ and $H=H(I(D)_{\ev{2}})$ both are co-chordal and $D_i$ are not induced subgraphs of $D$. Let $n>2$ and $D$ be non-empty. Since $G$ is co-chordal, then $G$ does not contain a cycle greater than $4$. Then the following three cases arise. \\
{\bf Case-1:} Suppose that $D$ contain at least one induced $4$-cycle and at least one induced $3$-cycle. Then either they are connected by a path, or they can share one vertex as in the graphs $G_2$ and $G_1$ in Figure \ref{fig8}. In both cases, we get a contradiction because induced subgraphs of the underlying graph $G$ are not co-chordal. Then case-$1$ never arises.\\
 {\bf Case-2:} Suppose that $D$ contain only induced $4$-cycles. Then $D$ be a bipartite graph. Then the required implication follows from Theorem \cite[Theorem 4.9]{kns25}.\\
{\bf Case-3:} Assume that $D$ contain only induced $3$-cycles, Then $D$ is a chordal graph. Therefore, by Theorem \ref{theorem-3.4}, we get $I(D)$ is vertex splittable.
 \end{proof}
 
\begin{lemma}\label{lemma-4.6}
   Let $D$ be a weighted oriented complete bipartite graph such that $H=H(I(D)_{\ev{2}})$ is co-chordal and $D_i$ are not induced subgraphs of $D$ for $i=1,\ldots,4$. Then
   $|V^+|\leq 1$.
\end{lemma}
 \begin{proof}
     Let $V(D)=A\cup B$ be a partition of $V(D)$ with $A=\{x_1,\ldots,x_m\}$ and $B=\{y_1,\ldots,y_n\}$.
     If possible, let $|V^+| \geq 2$. Let $x_a,x_b\in V^+$. Then both these vertices can be present in either the same partition or two different partitions. Suppose these are present in the same partition, say $A$, i.e., $x_a,x_b\in A\cap V^+$. Since $x_a,x_b$ not are source vertices this gives that $(y_i,x_a)\in E(D)$ and $(y_j,x_b)\in E(D)$ for some $i,j$. Thus, we should have edge  $(x_a,y_j)\in E(D)$ and  $(x_b,y_i)\in E(D)$ as $D_2$ is not an induced subgraph of $D$. Then $\{x_a,y_j\}$ and $\{x_b,y_i\}$ form an induced matching in $H$, which is a contradiction as $H$ is co-chordal.
     Now, assume that $x_a,x_b$ are present in different partitions. Without loss of generality, let $x_a\in A\cap V^+$ and $y_b\in B\cap V^+$. Since $x_a,y_b$ not are source vertices this implies that $(y_i,x_a)\in E(D)$ and $(x_j,y_b)\in E(D)$ for some $i,j$. Since $D$ is a complete bipartite graph, then either $(x_a,y_b)\in E(D)$ or $(y_b,x_a)\in E(D)$, but in both cases we get $D_1$ as an induced subgraph of $D$, which is a contradiction. Therefore $|V^+|\leq 1$.  
 \end{proof}
 
 \begin{lemma}\label{verpartite}
     Let $D$ be a weighted oriented complete $r$-partite graph with $r\geq2$ such that $H=H(I(D)_{\ev{2}})$ is co-chordal and $D_i$ are not induced subgraphs of $D$ for $i=1,\ldots,4$. Then $|V^+|\leq (r-1)$. Furthermore, the vertices of $V^+$ lie in the same partition of $V(D)$ or two different partitions of $V(D)$ with the distribution that one partition contains at most $(r-2)$ vertices and the other one contains at most one vertex of $V^+$. 
\end{lemma}
\begin{proof}
Let $V(D)=A_1\cup \dots \cup A_r$ be a partition of $V(D)$. We want to show that $|V^+|\leq (r-1)$. The lemma follows from Lemma \ref{lemma-4.6} if $r=2$. Assume now on wards $r\geq 3$. \\
If possible, let $|V^+|\geq r$. If the elements of $V^+$ are lying in three or more partitions of $V(D)$, then there exist $x,y,z \in V^+$ such that all these three vertices lie in different partitions. This implies that the induced subgraph $D[\{x,y,z\}]$ is either $D_3$ or $D_4$. This is a contradiction. Therefore the vertices of $V^+$ are lying in either in the same partition or two different partitions.\\ 
\noindent 
{\bf Case 1:} Suppose the vertices in $V^+$ are lying in the same partition $A_i$ of $V(D)$. Since $|A_i| \geq |V^+|\geq r$ and $D$ is an $r$-partition graph, then there exists a partition $A_j$ with $j\neq i$ such that the induced complete bipartite subgraph $D[A_i\cup A_j]$ has at least two elements of $V^+$. This is a contradiction by Lemma \ref{lemma-4.6}. 

\noindent 
{\bf Case-2:} If possible, let the vertices of $V^+$ are lying in two different partitions $A_i$ and $A_j$ of $V(D)$ for some $i\neq j$. 
Now, we want to show that at least one of these two partitions contains exactly one element of $V^+$. Suppose not. That is, each of these two partitions contains more than one vertex of $V^+$. Let $x_m,x_{m'}\in A_i\cap V^+$ and $x_n,x_{n'}\in A_j\cap V^+$. This implies that the induced complete bipartite subgraph $D[A_i\cup A_j]$ has at least two elements  of $\{x_m,x_{m'}, x_n,x_{n'}\}$. This is a contradiction by Lemma \ref{lemma-4.6}. Therefore, one partition, say $A_i$ contains only one vertex of $V^+$, and hence the other partition, say $A_j$ contains at least $r-1$ vertices of $V^+$. Let $x_{m} \in A_i \cap V^+$ and $x_{n} \in A_j \cap V^+$. Since $x_{m},x_{n}$ are non-source vertices, therefore $N^{-}_D(x_{m}) \neq \emptyset$ and $N^{-}_D(x_{n}) \neq \emptyset$. Let $z\in N^{-}_D(x_{m})$ and $u \in N^{-}_D(x_{n})$. If possible, $z,u$ lie in same partition, say, $z,u\in A_1$.  We have $(z,x_{m})\in E(D)$ and $(u,x_{n})\in E(D)$. This implies that $\{x_{m},u\}, \{x_{n},z\} \in E(H)$, otherwise $D_4$ is present in $D$ as an induced subgraph which is absurd. Then $\{x_{m}, u\}$ and $\{x_{n},z\}$ form an induced matching in $H$ because $(z,x_{m})\in E(D)$, $(u,x_{n})\in E(D)$ which leads to a contradiction if $z\neq u$ because $H$ is co-chordal. If $z=u$, then we get $D_4$ as an induced subgraph of $D$ which is a contradiction. Therefore $z,u$ must lie in different partitions of $V(D)$. Since $|A_j\cap V^+| \geq r-1$, and $x_n\in A_j\cap V^+$, then each element in $(A_j\cap V^+)\setminus \{x_n\}$ is adjacent to some vertex from other partition $A_k$, with $k\neq j$. Thus we have 

\begin{enumerate}
    \item  $|(A_j\cap V^+)\setminus \{x_n\}| \geq r-2$, 
    \item the $(r-2)$ number of elements of $(A_j\cap V^+)\setminus \{x_n\}$ must be adjacent to elements of $(r-1)$ number of disjoint partitions $A_1,\ldots,A_{j-1}, A_{j+1}, \ldots, A_r$, and 
    \item  $z\neq u$ and $z, u$ belong to two different partitions.  
\end{enumerate}
This implies that there exists a partition $A_l$ with $l\neq j$ such that the induced subgraph $D[A_l\cup A_j]$ is complete bipartite  and $|V^+(D[A_l\cup A_j])| \geq 2$. This is a contradiction by Lemma \ref{lemma-4.6}. 
 Thus in both cases, we get a contradiction, therefore $|V^+|\leq r-1$.\par

 \noindent 
 Hence the vertices of $V^+$ lie in the same partition or two different partitions, with the distribution that one partition contains at most $(r-2)$ vertices and the other one contains only one vertex of $V^+$. 
\end{proof}

\begin{remark}\label{partite}
    Suppose $D$ is a weighted oriented $r$-partite graph such that $H$ is co-chordal and $D_i$ are not induced subgraphs of $D$ for $1\leq i\leq 4$. Let $x,y\in V^+$ then $N^{-}_D(x) \neq \emptyset$ and $N^{-}_D(y) \neq \emptyset$. Let $z\in N^{-}_D(x)\cap A_i$ and $u \in N^{-}_D(y)\cap A_j$. Then $z,u$ lie in two different partitions. This implies that $z'\in N_D^+(V^+)$ for all $z'\in A_i$.
\end{remark}
Below we give a characterization for componentwise linearity of edge ideals of complete $r$-partite graphs.

\begin{theorem}\label{compartite}
     Let $D=(V(D), E(D), w)$ be a weighted oriented complete $r$-partite graph and $G$ its underlying simple graph. Then, the following statements are equivalent.
     \begin{enumerate}
         \item $I(D)$ is componentwise linear;
         \item $I(D)$ is vertex splittable;
         \item $I(D)$ has linear quotient property;
         \item both $G$ and $H=H(I(D)_{\ev{2}})$ are co-chordal and $D_i$ are not induced subgraphs of $D$ for $i=1,\ldots,4$.
     \end{enumerate}
\end{theorem}     
\begin{proof}
$(2) \implies (3) \implies (1)$ is easy. \\ 
\noindent $(1)\implies (4)$: Let $I(D)$ be componentwise linear. Then by Theorem \ref{corollary-2.7}, we get $D_i$ are not induced subgraphs of $D$ for $i\in \{1,\ldots,4\}$. Also using Lemma \ref{thmcochordal}, we obtained that $G$ is co-chordal. Since $I(D)$ is componentwise linear then $H=H(I(D)_{\ev{2}})$ has linear resolution. Therefore, by { Fr\"{o}berg's Theorem} we get $H$ is co-chordal.\\
\noindent $(4)\implies (2)$: We proceed by induction on the number of vertices $n$ of $D$. Let $n\geq 3$ and $D$ be non-empty. Note that  $N_G(x)$ is a minimal vertex cover of $G$ for any $x \in G$. By Lemma \ref{verpartite}, we have at most $r-1$ vertices of $D$ which are in $V^+$ and we get a vertex $x\in V(D)$  with $w(x)=1$. Then we can write as follows $I(D)=xI_1+I(D\setminus x)$,
   where $$I_1=\left(z^{w(z)}\mid z\in N_{D}^{+}(x)\cap V^{+}\right)+ (z\mid z\in  \left (N_{D}^{+}(x)\setminus V^{+}) \cup N_{D}^{-}(x)\right).$$ 
   Using Lemma \ref{lemma-3.1}, we have $|N_D^+(x)\cap V^+|\leq 1$. Then by Lemma \ref{lemma-2.8} , we get that $I_1$ is vertex splittable. Also, $I(D\setminus x)$ is vertex splittable by the induction hypothesis. If $|N_D^+(x)\cap V^+|=\emptyset$ then $I(D\setminus x)\subseteq I_1$.  Hence, $I(D)$ is vertex splittable. If $|N_D^+(x)\cap V^+|= 1$, without loss of generality assume $x_1\in N_D^+(x)\cap V^+$. If $N_D^+(x_1)\setminus N_D(x) = \emptyset$ then $I(D\setminus x)\subseteq I_1$. Hence, $I(D)$ is vertex splittable, as required.\par
 Suppose $x_1$ is not a sink vertex i.e., $N_D^+(x_1)\setminus N_D(x) \neq \emptyset$. Let $z\in N_D^+(x_1)\setminus N_D(x)$. This implies that, $z$ and $x$ are in the same partition and $w(z)=1$ as $D_1$ cannot be an induced subgraph of $D$. Since $(x,x_1)\in E(D)$ and $z,~x$ are in the same partition  then by Remark \ref{partite} we have $z\in N_D^+(V^+)$. Therefore, we can write $I(D)=zN_D(z)+I(D\setminus z)$.  Since $N_{D}(z)$ is a minimal vertex cover of $D$ implies $I(D\setminus z)\subseteq (y~|~ y\in N_D(z))$. Also, $I(D\setminus z)$ is vertex splittable by the induction hypothesis. Hence, $I(D)$ is vertex splittable.
\end{proof}

\section{Powers of componentwise linear house-free weighted oriented graphs}\label{secpowervertexsplit}

 In this section, we study the componentwise linearity of powers of edge ideals of house-free weighted oriented graphs. 

 \noindent 
Recall that if $I=(f_1,\ldots,f_r)$ is a monomial ideal in $R$ and $g\in R$ a monomial, then $(I:g)$ is generated by the monomials $f_i:g = \frac{f_i}{\text{gcd}(f_i,g)}$, for $1\leq i\leq r$. 
 
 In general, if $I(D)$ is componentwise linear, then $I(D)^2$ need not be componentwise linear. See the below example. 
 
\begin{example}\label{power}
    Let $I(D)=(x_1x_2,x_2x_3,x_3x_1,x_1x_4^2,x_2x_5^2,x_3x_6^2)$. Then $I(D)$ is componentwise linear, but $I(D)^2$ is not componentwise linear. By Macaulay2 \cite{gs}, one can check that $\reg(I(D)^k)=3k+1$, for all $k\geq 2$. This implies that $I(D)^k$ is not componentwise linear, for all $k\geq 2$. 
\end{example}

\begin{remark}\label{remark-4.1}
   Let $D$ be a weighted-oriented graph. Suppose $I(D)=x I_1 + I_2,$ is a vertex splittiable,  where
\(
I_1 = \left( \{x_i^{w(x_i)} \mid x_i \in \N_D^+(x)\} \right)+ \left ( z~|~z\in \N_D^-(x) \right)\), \(I_2 = I(D \setminus x),\) with $I_2 \subseteq I_1$, and $N_D(x)$ is a minimal vertex cover. By Lemma~\ref{lemma-3.1}, we have $|\N_D^+(x) \cap V^+| \leq 1$, so $I_1$ is of the form
\(
I_1 = (x_1, \ldots, x_i^{w(x_i)}, \ldots, x_n, y_1, \ldots, y_r),
\)
where $w(x_i) \geq 1$ for some $i$, and $y_j \in \N_D^-(x)$. Then,
\begin{align*}
I(D)^2 &= (xI_1 + I_2)^2  
       = x I_1 I(D) + I_2^2.
\end{align*}
Therefore, by \cite[Lemma~3.2]{fmr24}, if both $I_1 I(D)$ and $I_2^2$ have linear quotients, then $I(D)^2$ also has linear quotients.
\end{remark}
\begin{proposition}\label{D_4}
Let $D_4$ is the weighted oriented graph as in Figure \ref{fig1} with the edge ideal $I(D_4)=(x_1x_2^{w_2},x_1x_3^{w_3},x_3x_2^{w_2})$ where $w_2\geq2$ and $w_3\geq2$. Then, 
\begin{equation*}
 \reg(I(D_4)^k)=(k-1)(w+1)+w_2+w_3, \text{ for all } k\geq 1,  
\end{equation*}
where $w=\text{max}\{w_2,w_3\}$. In particular, $I(D_4)^k$ is not componentwise linear for any $k\geq 1$. 
\end{proposition}
\begin{proof}
    We prove by induction on $k$. If $k=1$, then it is easy to see $\reg(I(D_4))=w_2+w_3$. Now, assume that $k\geq 2$. We have 
\begin{align*}
    I(D_4)^k&=(x_1x_2^{w_2},x_1x_3^{w_3},x_3x_2^{w_2})^k
   =\displaystyle\sum_{i=0}^k(x_1x_3^{w_3})^{k-i}(x_1x_2^{w_2},x_3x_2^{w_2})^i \\
   &=(x_1x_3^{w_3})^k+\displaystyle\sum_{i=1}^k(x_1x_3^{w_3})^{k-i}(x_1x_2^{w_2},x_3x_2^{w_2})^i \\
   & =(x_1x_3^{w_3})^k+(x_1x_2^{w_2},x_3x_2^{w_2})I(D)^{k-1} = K+J \text{ (say)}
\end{align*}
 where $K=(x_1^kx_3^{kw_3})^k$ and $J=(x_1x_2^{w_2},x_3x_2^{w_2})I(D)^{k-1}$. Let $J^{\PP},~K^{\PP}$ and $(I(D_4)^k)^{\PP}$ denote the polarizations of $J,~K$ and $I(D_4)^k$ respectively. Thus we have $(I(D_4)^k)^{\PP}=J^{\PP}+K^{\PP}$. One can check that  
    \begin{equation*}
     J^{\PP}\cap K^{\PP} =\Bigg (\prod_{j=1}^k x_{1j}\prod_{j=1}^{w_2}x_{2j}\prod_{j=1}^{kw_3}x_{3j} \Bigg ).
    \end{equation*}
    Since $K^{\PP}$ has linear resolution and the variable $x_{3,k{w_3}}$ in $K^{\PP}$ 
cannot divide any element in $\mathcal{G}(J^{\PP})$, then we get that 
$(I(D_4)^k)^{\PP} = J^{\PP} + K^{\PP}$ is a Betti splitting by \cite [Lemma 2.2]{x21}.
Therefore,
\begin{align}\label{eqn1}
    \reg (I(D_4)^k)=\reg(I(D_4)^k)^{\PP}=\text{max}\{\reg (K^{\PP}),~\reg (J^{\PP}),~\reg (K^{\PP}\cap K^{\PP})-1\}. 
\end{align}
Now, consider the short exact sequence 
\begin{equation*}
 0 \rightarrow x_1x_3x_2^{w_2}I(D_4)^{k-1} \rightarrow x_1x_2^{w_2}I(D_4)^{k-1}\oplus  x_3x_2^{w_2}I(D_4)^{k-1} \rightarrow J \rightarrow 0
\end{equation*}
Let, $A=x_1x_3x_2^{w_2}I(D_4)^{k-1}$, $B= x_1x_2^{w_2}I(D_4)^{k-1}\oplus  x_3x_2^{w_2}I(D_4)^{k-1}$. Then $\reg (B)=\reg (I(D_4))^{k-1}+w_2+1$ and $\reg (A)=\reg (I(D_4)^{k-1})+w_2+2$, hence $\reg (A)=\reg (B)+1$. Then by the regularity lemma \cite[Corollary 20.19]{e95} applied to the above short exact sequence, we get $\reg J= \reg (A)-1$. Thus $ \reg J=\reg (I(D_4)^{k-1})+w_2+1$.\\
If $w_2\geq w_3$ then by induction we have $\reg (I(D_4)^{k-1})=(k-2)(w_2+1)+w_2+w_3$. Then $\reg J=(k-2)(w_2+1)+w_2+w_3+w_2+1=(k-1)(w_2+1)+w_2+w_3$. Since $\reg J=\reg J^{\PP}$ then $\reg J^{\PP}=(k-1)(w_2+1)+w_2+w_3$. Therefore using Equation \eqref{eqn1} we have,
\begin{align*}
   \reg (I(D_4)^k) &=\text{max}\{k(w_3+1), (k-1)(w_2+1)+w_2+w_3, k(w_3+1)+w_2-1\}\\
   &=\text{max}\{(k-1)(w_2+1)+w_2+w_3, k(w_3+1)+w_2-1\}\\
   &=(k-1)(w_2+1)+w_2+w_3.
\end{align*}
If $w_3> w_2$, then similarly we get, 
  \begin{align*}
   \reg (I(D_4)^k)&=\text{max}\{k(w_3+1), (k-1)(w_3+1)+2w_2, k(w_3+1)+w_2-1\}\\
  &=k(w_3+1)+w_2-1=(k-1)(w_3+1)+w_3+w_2. 
\end{align*}
Therefore $\reg(I(D_4)^k)=(k-1)(w+1)+w_2+w_3, \text{ for all } k\geq 1.$ This implies that $\reg(I(D_4)^k) > k(w+1)$, for all $k\geq 1$. This gives that $I(D_4)^k$ is not componentwise linear for any $k\geq 1$. 
\end{proof}

\begin{theorem}\label{D_i}
   Let $D$ be a weighted oriented graph. If $I(D)^k$ is componentwise linear for some $k\geq1$, then $D_i$ are not induced subgraphs of $D$ for all $i\in \{1,\ldots,4\}$. 
\end{theorem}
\begin{proof}
    Let $I(D)^k$ is componentwise linear for some $k\geq1$. Suppose $D_i$ is an induced subgraph of $D$ for some $1\leq i \leq 4$. Then by Lemma \ref{inducedpower}, we get $I(D_i)^k$ also componentwise linear. We will show the required result in four individual cases: \\
{\bf Case-1:} Suppose $D_i=D_1$. Assume the notation for $D_1$ as in Figure \ref{fig1}. Since $w(x_2)>1$ and $w(x_3)>1$, by \cite[Lemma 3.4]{x21}, we get,
\begin{equation*}
 \reg(I(D_1)^k)=\displaystyle\sum_{i=1}^3w(x_i)-1+(k-1)(w+1),   
\end{equation*}
where $w=\text{max}\{w_2,w_3\}$. Since $w_2>1$ and $w_3>1$, we have $\reg(I(D_1)^k)>k(w+1)$.
Therefore, $I(D_1)^k$ is not componentwise linear. Thus $D_i \neq D_1$. \\
{\bf Case-2:} Suppose $D_i=D_2$. Assume the notation for $D_2$ as in Figure \ref{fig1}. Then  by \cite[Theorem 3.3]{x21} we get,
 \begin{align*}
 \reg(I(D_2)^k)&=\displaystyle\sum_{i=1}^3w(x_i)-|E|+1+(k-1)(w+1)\\
 &=w_2+w_3+(k-1)(w+1),
 \end{align*}
 where $w=\text{max}\{w_2,w_3\}$. Thus, by $w_2>1$ and $w_3>1$, we have $\reg(I(D_2)^k)=w_2+w_3+(k-1)(w+1)>k(w+1)$. Therefore, $I(D_2)^k$ is not componentwise linear. Thus $D_i \neq D_2$.\\
 {\bf Case-3:} Suppose $D_i=D_3$. Assume the notation for $D_3$ as in Figure \ref{fig1}. Then from \cite[Theorem 4.5]{wg23}, it follows that
 \begin{align*}
 \reg(I(D_3)^k)&=\displaystyle\sum_{i=1}^3w(x_i)-|E|+1+(k-1)(w+1)\\
 &=w_2+w_3+(k-1)(w+1),
 \end{align*} 
  where $w=\text{max}\{w_2,w_3\}$.
 Thus, by $w_2>1$ and $w_3>1$, we have $\reg(I(D_2)^k)=w_2+w_3+(k-1)(w+1)>k(w+1)$. Therefore, $I(D_3)^k$ is not componentwise linear. Thus $D_i \neq D_3$.\\
 {\bf Case-4:}  Suppose $D_i=D_4$. Then by Proposition \ref{D_4}, we get $I(D_4)^k$ is not componentwise linear. Thus $D_i \neq D_4$.\\
 Hence $D_1,D_2,D_3,D_4$ cannot be induced subgraphs of $D$.
\end{proof}

Below we show that if $D$ is house-free and $I(D)$ is not componentwise linear, then $I(D)^k$ is not componentwise linear, for all $k\geq 2$. 

\begin{corollary} \label{cor1}
 Let $D$ be a weighted oriented graph and $G$ its underlying simple graph such that $G$ is $G^{\prime}$-free (that is, house-free),  where $G'$ is as in Figure \ref{fig8} such that $G$ and $H=H(I(D)_{\ev{2}})$ are co-chordal. If $I(D)^k$ is componentwise linear for some $k\geq 2$, then $I(D)$ is componentwise linear.    
\end{corollary}
 \begin{proof}
    Let $D$ be a weighted oriented house-free graph such that the underlying graph $G$ and $H=H(I(D)_{\ev{2}})$ are co-chordal.  If $I(D)^k$ is componentwise linear for some $k\geq 2$, then by Theorem \ref{D_i} we get $D_i$ are not induced subgraphs of $D$ for all $1\leq i \leq 4$. Therefore by Theorem \ref{thm3} we get that $I(D)$ is componentwise linear.
 \end{proof}   

\noindent 
 The converse of Corollary \ref{cor1} is not true, see Example \ref{power}.

\section{Powers of componentwise linear weighted oriented complete $r$-partite graphs} \label{sec5}

In this section, we study the componentwise linearity of powers of edge ideals of complete $r$-partite weighted oriented graphs. We show that if $D$ is a complete $r$-partite graph, then $I(D)$ is not componentwise linear if and only if $I(D)^k$ is not componentwise linear, for all $k\geq 2$. 

\begin{lemma}\label{lemma-4.4}
Let $D$ be a weighted oriented graph and $I = I(D)$ its edge ideal and let $P \subset S$ be a monomial prime ideal. Suppose $I$ has linear quotients. Then $PI$ has linear quotients.
\end{lemma}
\begin{proof}
Without of loss of generality, we may assume $P = (x_1, \ldots, x_t)$. Let $I=(u_1, \ldots, u_m)$, where $u_1, \ldots, u_m$ is a linear quotient order of $I$. We prove the lemma by induction on $m$. If $m = 1$, then it is easy to see that $x_1 u_1, \ldots, x_t u_1$ is a linear quotients order for $PI$, as required. Assume $m \geq 2$. Let $L = (u_1, \ldots, u_{m-1})$. Then $I = L +(u_m)$, and $L$ is the edge ideal of a weighted oriented graph with linear quotients. Therefore by induction hypothesis, $PL$ has a linear quotient order, say, $v_1, \ldots, v_h$ and $PL=(v_1, \ldots, v_h)$. Let $\mathcal{G}(PI) \setminus \mathcal{G}(PL) =\{x_{j_1}u_m, \ldots, x_{j_s}u_m \}$, for $1 \leq j_1 < \ldots < j_s \leq t$. We claim that
$v_1, \ldots, v_h, x_{j_1}u_m, \ldots, x_{j_s}u_m$ is a linear quotient order for $PI$. Note that $PI=(v_1, \ldots, v_h)+(x_{j_1}u_m, \ldots, x_{j_s}u_m )$. Since by induction, $v_1, \ldots, v_h$ is a linear quotients order of $PL$, it remains to show that $((v_1, \ldots, v_h, x_{j_1}u_m, \ldots, x_{j_{i-1}}u_m) : x_{j_i}u_m)$ is generated by variables, for all $1 \leq i \leq s$. It is clear that $x_{j_p}u_m : x_{j_i}u_m = x_{j_p}$, for all $1 \leq p < i$. Let $v_{\ell} = x_p u_q$, for some $1 \leq p \leq t$ and for some $1 \leq q < m$. Consider the monomial $v_{\ell} : x_{j_i}u_m$,  if its degree is $1$, then there is nothing to prove. Hence assume that $\deg(v_{\ell} : x_{j_i}u_m)\geq 2$. Let $u_q = x_r x_s^{w(x_s)}$. Then at least one of the variables $x_r$ and $x_s^{w(x_s)}$ divides $v_{\ell} : x_{j_i}u_m$, this imply that either $x_r$ or $x_s$ divides $v_{\ell} : x_{j_i}u_m$ because if both $x_r$ and $x_s$ does not divide  $v_{\ell} : x_{j_i}u_m$, then $\deg(v_{\ell} : x_{j_i}u_m) = 1$.\\
\textbf{Case-1:} Suppose $x_r$  divide $v_{\ell} : x_{j_i}u_m$. Consider the monomial $u_q : u_m$ and $I$ has linear quotients, then there exists $k < m$ such that $u_k : u_m$ is equal to a variable and that variable divides the monomial $u_q : u_m$. Thus $u_k : u_m$ divides $x_r x_s^{w(x_s)}$. If $u_k : u_m = x_r$, then $x_{j_i}u_k : x_{j_i}u_m = x_r$. Notice that $x_{j_i}u_k \in \mathcal{G}(PL)$. So in this case we are done. Let us assume $u_k : u_m = x_s$, if $x_s$ divides $(v_{\ell} : x_{j_i}u_m)$ then the same argument as before can be applied. If $x_s$ does not divide $v_{\ell} : x_{j_i}u_m$, then $v_{\ell} : x_{j_i}u_m = x_p x_r$. Since $u_k : u_m = x_s$, then $x_s$ does not divide $u_m$. These imply that $j_i = s$ and hence $u_k$ divides $x_{j_i}u_m$. Therefore, $x_p u_k : x_{j_i}u_m = x_p$ divides $v_{\ell} : x_{j_i}u_m$ and $x_p u_k \in \mathcal{G}(PL)$.\\
\textbf{Case-2:} Assume $x_s$ divide $v_{\ell} : x_{j_i}u_m$. Consider the monomial $u_q : u_m$ and $I$ has linear quotients, then there exists $k < m$ such that $u_k : u_m$ is a variable that divides $u_q : u_m$, and so $u_k : u_m$ divides $x_r x_s^{w(x_s)}$. If $u_k : u_m = x_s$, then $x_{j_i}u_k : x_{j_i}u_m = x_s$. Notice that $x_{j_i}u_k \in \mathcal{G}(PL)$. So in this case we are done. Let us assume $u_k : u_m = x_r$ if $x_r$ divides $(v_{\ell} : x_{j_i}u_m)$, then the same argument as before can be applied. If $x_r$ does not divide $v_{\ell} : x_{j_i}u_m$, then $v_{\ell} : x_{j_i}u_m = x_p x_s^{w(x_s)}$. Since $u_k : u_m = x_s$, and $x_s$ does not divide $u_m$, this imply that $j_i = s$ and hence $u_k$ divides $x_{j_i}u_m$. Therefore, $x_p u_k : x_{j_i}u_m = x_p$ divides $v_{\ell} : x_{j_i}u_m$ and $x_p u_k \in \mathcal{G}(PL)$. This completes the proof. 
\end{proof}

\noindent 
Lemma \ref{lemma-4.4} need not be true for any ideal $I$. See the example below.  
\begin{example}
Let $S=K[a,b,c,d]$, where $K$ is a field. Let $I=(a^2b,abc,bcd,cd^2)$, and $P=(a,b)$. Then one can check that $I$ has linear quotients, but $PI$ does not have linear resolution, and hence $PI$ does not have linear quotients. 
\end{example}

\begin{lemma}\label{lemma-4.5}
  Let $D$ be a weighted oriented graph and $I(D)$ its ideal. Suppose  $I(D)=x(x_1,\ldots,x_n)+I(D\setminus x)$, where $\N_D(x)=\{x_1,\ldots,x_n\}$ is minimal vertex cover and $w(x)=1$. Then $I(D)^2$ has the linear quotient property if $I(D)$ has linear quotient property. 
\end{lemma}
\begin{proof} Under the assumptions of the lemma, one can get that  $I(D)^2=x(x_1,\ldots,x_n)I(D)+I(D\setminus x)^2$. To show that $I(D)^2$ has linear quotient property, we prove by induction on $k=\vert V(D) \vert$. If $k=2$, then it is easy to see that $I(D)^2$ has linear quotient property. Assume  $k \geq 2$. By induction hypothesis, $I(D\setminus x)^2$  has linear quotients. By Lemma \ref{lemma-4.4}, we have that $(x_1,\ldots,x_n)I(D)$ has linear quotient property. Therefore by Remark \ref{remark-4.1}, we have that $I(D)^2$ has linear quotient property.
\end{proof}

\begin{lemma}\label{comquotient1}
  Let $D$ be a weighted oriented complete $r$-partite graph with the underlying simple graph $G$. If $I(D)$ has linear quotients then $I(D)$ has linear quotient order $$u_1,\ldots,u_t,u_{t+1},\ldots,u_h,$$ of generator of $I(D)$, such that $\deg(u_1)=\ldots=\deg(u_t)=2$ and $\deg(u_i)\geq 3$, for all $t+1\leq i\leq h$ with $u_1,\ldots,u_t$, the linear quotient order of $I(D)_{\ev{2}}$ and $u_{t+1},\ldots,u_h$ can be permuted.
  \end{lemma}
  \begin{proof}
    Let $D$ be a weighted oriented complete $r$-partite graph with the partition $V(D)=A_1\cup\ldots \cup A_r$. Assume $I(D)$ has linear quotient property. By Lemma \ref{increasing-admissible} we have a degree increasing linear quotient order $u_1,\ldots,u_h$ (say) of $\mathcal{G}(I(D))$ such that $\deg(u_1)\leq \ldots\leq \deg(u_t)\leq \deg(u_{t+1})\leq \ldots\leq \deg(u_h)$. Let  $\deg (u_i)=2$ for $1\leq i\leq t$. Then we have $I(D)_{\ev{2}}$ has linear quotient order $u_1,\ldots,u_t$ . Since $I(D)$ has linear quotient property, then by Theorem \ref{corollary-2.7} we get $D_i$ are not induced subgraphs of $D$ for all $i\in \{1,\ldots,4\}$. Then using Lemma \ref{verpartite}, we get that $|V^+|\leq r-1$ and the vertices of $V^+$ lie in at most two partitions of $V(D)$. Denote $V^+=\{x_1,\ldots,x_\alpha\}$. Note that $\alpha \leq r-1$. Therefore, we have the following two cases\\
\textnormal{(i)} The vertices of $V^+$ lie in  two different partitions.\\
\textnormal{(ii)} The vertices of $V^+$  lie in same partition.\\
\textbf{\underline{Proof of (i):}} Assume the vertices of $ V^+$ lie in the two different partition. Assume  $x_1,\ldots,x_{\alpha-1}\in A_i\cap V^+$ and $x_\alpha\in A_j\cap V^+$ for some $i\neq j$. Let $\mathcal{G}(I(D))=S \cup S_{1} \cup \ldots \cup S_\alpha$, where
\begin{align*}
 S=&\{xy:xy\in \mathcal{G}(I(D))\},~ 
    S_{k}=\{xx_{k}^{w_{k}}\in \mathcal{G}(I(D)):x\in N_D^-(x_{k})\}, ~\text{for}~~ 1\leq k\leq \alpha.
  \end{align*}
 Note that if we assume the notation as in the statement, we have $S_1=\{u_1,\ldots,u_t\}$ and $S_{1}\cup \ldots\cup S_\alpha= \{u_{t+1}, \dots, u_h\}$, where $u_1,\ldots,u_t$ is a linear quotient order for $I(D)_{\langle 2\rangle}$.  
We want to show that $u_1,\ldots,u_t,u_{t+1},\ldots,u_h$ is a linear quotient order of $\mathcal{G}(I(D))$. In fact, we show that if $u_m, u_n \in S_k$, for any $1\leq k \leq \alpha$, there exists $u_l\in S$ such that $u_l:u_n$ is a variable and this variable divides $u_m:u_n$. This we are proving in various cases. \\
  \textbf{Case-1:} Let $u_m\in S$ precedes $u_n\in S_{k}$ ($1\leq k\leq \alpha-1$), say $u_m=xy$ and $u_n=zx_{k}^{w_{k}}$. If $u_m:u_n$ is a variable, then we are done. Suppose $u_m:u_n$ is not a variable. Then $x,y\neq \{z,x_k\}$. Assume $z=x_\alpha\in A_j\cap V^+$. Since $xy$ is an edge of $D$ then one of $x,y\notin A_j$, without loss of generality $x\notin A_j$. Thus, either $x\in N_D^-(x_\alpha)$ or  $x\in N_D^+(x_\alpha)$. This implies that either $(x,x_\alpha)\in E(D)$ or $(x_\alpha,x)\in E(D)$. Let $(x,x_\alpha)\in E(D)$ then $(x_k,x)\in E(D)$ as we have $(x_\alpha,x_{k})\in E(D)$ and $D_1$ is not an induced subgraph of $D$. This implies that $x_{k}x\in S_1$, say $u_l=x_{k}x$. Then $u_l:u_n=x$ divides $u_m:u_n$, as required. Let $(x_\alpha,x)\in E(D)$. Since $x_{k}\in N_D^+(x_\alpha)\cap V^+ $ then using Lemma \ref{lemma-3.1}, we have $xx_\alpha\in S_1$, say $u_l=xx_\alpha$. Then $u_l:u_n=x$ divides $u_m:u_n$, as required. Assume  $z\neq x_\alpha$. Then $w(z)=1$ and $z\notin A_i$. Let $z\in A_l$ for some $l\neq i$. We have one of $x$ or $y$ not in $A_l$ as $xy$ is an edge of $D$, without loss of generality assume $x\notin A_l$. If $w(x)>1$, then using Lemma \ref{lemma-3.1} we get $(x,z)\in E(D)$ because $x_{k}\in N_D^+(z)\cap V^+$. This implies that $xz\in S_1$. If $w(x)=1$, then it is easy to see that $xz\in S_1$ as we have $w(z)=1$. Thus in both cases, we get $xz\in S_1$, say $u_l=xz$. Then $u_l:u_n=x$ divides $u_m:u_n$, as required. \par
\noindent
\textbf{Case-2:} Let $u_m\in S_{k'}$ precedes $u_n\in S_{k}$ with $u_m=xx_{k'}^{w_{k'}}$ and $u_n=zx_{k}^{w_{k}}$ for $k,k'\in \{1,\ldots,\alpha-1\}$. Then $x\neq z$ as by Lemma \ref{lemma-3.1} we have $|N_D^+(x)\cap V^+|\leq 1$. Also using Remark \ref{partite}, we get that $x,z$ are lying in different partitions. Assume $u_m:u_n$ is not a variable. This is same as Case-1 by interchanging the role of $x_{k'}$ and $y$.\par
\noindent
\textbf{Case-3:} Let $u_m\in S_\alpha$ precedes $u_n\in S_{k}$ with $u_m=xx_\alpha^{w_\alpha}$ and $u_n=zx_{k}^{w_{k}}$ for $k\in \{1,\ldots,\alpha-1\}$. Then $x\neq z$ as by Lemma \ref{lemma-3.1} we have $|N_D^+(x)\cap V^+|\leq 1$. Also using Remark \ref{partite}, we get that $x,z$ are lying in different partitions. Assume $u_m:u_n$ is not a variable. This is same as Case-1 by interchanging the role of $x_\alpha$ and $y$.\par
\noindent
\textbf{Case-4:} Let $u_m\in S_{k}$ precedes $u_n\in S_\alpha$ with $u_m=xx_{k}^{w_{k}}$ and $u_n=zx_\alpha^{w_\alpha}$ where $1\leq k\leq \alpha-1$. Then $x\neq z$ because by Lemma \ref{lemma-3.1} we have $|N_D^+(x)\cap V^+|\leq 1$. Let $z\in A_l$ for some $l\neq j$ then by Remark \ref{partite} $x\notin A_l$. Let $z\in A_i\cap V^+$ then, either $x\in N_D^-(z)$ or  $x\in N_D^+(z)$ this implies that either $(x,z)\in E(D)$ or $(z,x)\in E(D)$. Let $(x,z)\in E(D)$ then $(x_\alpha,x)\in E(D)$ as we have $(z,x_\alpha)\in E(D)$ and $D_1$ not an induced subgraph of $D$. This implies that $x_\alpha x\in S_1$, say $u_l=x_\alpha x$. Then $u_l:u_n=x$ divides $u_m:u_n$, as required. Let $(z,x)\in E(D)$. Since $x_\alpha\in N_D^+(z)\cap V^+ $ then using Lemma \ref{lemma-3.1}, we have $xz\in S_1$, say $u_l=xz$. Then $u_l:u_n=x$ divides $u_m:u_n$, as required. Assume  $z\notin A_i\cap V^+$ then $w(z)=1$. If $w(x)>1$ then using Lemma \ref{lemma-3.1} we get $(x,z)\in E(D)$ because $x_\alpha\in N_D^+(z)\cap V^+$. This implies that $xz\in S_1$. If $w(x)=1$ then $xz\in S_1$ as $w(z)=1$. In both cases, we get $xz\in S_1$, say $u_l=xz$. Then $u_l:u_n=x$ divides $u_m:u_n$, as required.\par
\noindent
\textbf{Case-5:} Let $u_m\in S$ precedes $u_n\in S_\alpha$, say $u_m=xy$ and $u_n=zx_{\alpha}^{w_{\alpha}}$. If $u_m:u_n$ is a variable, then we are done. Suppose $u_m:u_n$ is not a variable. Then $x,y\neq \{z,x_\alpha\}$. Let $z\in A_l$ for some $l\neq j$. Since one of $x,y\notin A_l$ as $xy$ is an edge of $D$, without loss of generality $x\notin A_l$.  This is same as Case-4 by interchanging the role of $x_k$ and $y$.\par
\noindent
\textbf{\underline{Proof of (ii):}} In this case, we get only Case-1 and Case-2 of proof (i). Proceeding as like in Proof of (i), we get the required result.
\end{proof}  

 \begin{theorem}\label{theorem-4.9}
 Let $D$ be a weighted oriented complete $r$-partite graph with $r\geq2$, and let \( I(D) \) denote its edge ideal. Then \( I(D)^2 \) has linear quotient property if  \( I(D) \) has linear quotient property.    
\end{theorem} 
\begin{proof}
   Let $D$ be a weighted oriented complete $r$-partite graph with the partition $A_1=\{x_{11}\ldots x_{1n_1}\}$,\ldots,$A_r=\{x_{r1},\ldots,x_{rn_r}\}$ and $I(D)$ has linear quotient property. Let $I(D)=(u_1,\ldots,u_h)$, where $u_1,\ldots,u_h$ is a linear quotient order as in Lemma \ref{comquotient1}. Now, by Lemma \ref{verpartite}, we have $|V^+|\leq r-1$ and the vertices of $V^+$ are lying in at most two different partitions. Denote $|V^+|=\alpha$. We have $\alpha \leq r-1$. Assume $x_{1},\ldots,x_{\alpha}\in V^+$. Using Lemma \ref{verpartite}, we get a vertex $x$ with $w(x)=1$ such that $(x,x_k)\in E(D)$ for some $1\leq k \leq \alpha$. Without loss of generality, let $x=x_{11}\in A_1$ and $x_k=x_{21}\in A_2$ then we get $(x,x_{21})\in E(D)$. Since $I(D)$ has linear quotient property, then by Theorem \ref{compartite}, $I(D)$ is vertex splittable and if $N_D^+(x_{21})\setminus N_D(x_{11})=\emptyset$, then we can write $I(D)=x_{11}(x_{21}^{w(x_{21})},x_{22}, \ldots,x_{2n_2},\ldots, x_{rn_r})+I(D\setminus x_{11})$. Thus, we have 
\begin{equation} \tag{*}\label{eq}
 I(D)^2=x_{11}PI(D)+I(D\setminus x_{11})^2,~~\text{where}~~ P=(x_{21}^{w(x_{21})},\ldots,x_{2n_2},\ldots,x_{r1},\dots x_{rn_r}).
 \end{equation}
By Lemma \ref{comquotient1}, we can take $u_h=x_{11}x_{21}^{w(x_{21})}$. To show $I(D)^2$ has linear quotient property, we prove by induction on $n=|V(D)|$. If $n=2$, then it is trivial. Assume $n\geq 3$. First we show that $PI(D)$ has linear quotients. We prove this by induction on $h$. If $h=1$, then it is easy to see that $x_{22}u_1,\ldots,u_1x_{rn_r},x_{21}^{w(x_{21})}u_1$ is a linear quotient order for $PI(D)$. Assume that $h\geq 2$. Let  $L=(u_1,\ldots,u_{h-1})$. By induction hypothesis, $PL$ has linear quotient property with linear quotient order, say, $v_1,\ldots,v_m$ and $PL=(v_1,\ldots,v_m)$. Consider 
\begin{align*}
Px_{11}x_{21}^w=&(x_{21}^{w(x_{21})},\ldots,x_{2n_2},\ldots,x_{r1},\dots x_{rn_r})x_{11}x_{21}^{w(x_{21})}\\
= &(x_{11}x_{22},\ldots,x_{11}x_{2n_2},\ldots ,x_{11}x_{r1},\dots x_{rn_r}x_{11})x_{21}^{w(x_{21})}+( x_{11}x_{21}^{2w(x_{21})}) \\
=& Ux_{21}^{w(x_{21})}+(x_{11}x_{21}^{2w(x_{21})}), 
\end{align*}
where $U=(x_{11}x_{22},\ldots,x_{11}x_{2n_2},\ldots ,x_{11}x_{r1},\dots x_{rn_r}x_{11})$. Since $D$ is a complete $r$-partite graph and $x_{21}\in N_D^+(x_{11})\cap V^+$ then by Lemma \ref{lemma-3.1} we have $U\subseteq L$. Then $Ux_{21}^w\in PL$ as $x_{21}^w\in P$ and $U\subseteq L$.
Therefore, the above expression for $PI(D)$ implies that $PI(D)=PL+(x_{11}x_{21}^{2w(x_{21})})$.
It remains to show that $((v_1,\ldots,v_m):x_{11}x_{21}^{2w(x_{21})})$ is generated by variables. Since the highest power of $x_{21}$ in any $v_l$ is at most $2w(x_{21})$ and the highest power of $x_{11}$ in any $v_l$ is at most $1$. Therefore, $x_{11}$ and $x_{21}$ does not divides $(v_l:x_{11}x_{21}^{2w(x_{21})})$, for all $l=1,\ldots,m $, i.e $(v_l:x_{11}x_{21}^{2w(x_{21})})\not\subset (x_{11},x_{21})$, for all $l=1,\ldots,m$.\\
 Since, either $x_{1i}\in N_D^+(x_{21})$ or  $x_{1i}\in N_D^-(x_{21})$ and $x_{21}^{w(x_{21})}\in P$, we get that $x_{1i}x_{21}x_{21}^{w(x_{21})}\in \mathcal{G}(PL)$ or $x_{1i}x_{21}^{2w(x_{21})}\in \mathcal{G}(PL)$ for all $2\leq i\leq n_1$. Also we have $x_{i{j_i}}x_{11}x_{21}^w\in \mathcal{G}(PL)$ for all $2\leq i\leq r$ and $1\leq {{j_i}}\leq n_i$ where $(i,j_i)\neq (2,1)$.  Therefore, these are some $ v_i$'s in the linear quotient order, which implies that 
 $$((v_1,\ldots,v_m):x_{11}x_{21}^{2w(x_{21})})=(\{x_{11},\ldots,x_{1n_1},\ldots,x_{r1},\ldots,x_{rn_r}\}\setminus \{x_{11},x_{21}\}).$$ 
 Therefore, $PI(D)$ has linear quotient property. By induction, we have $I(D\setminus x_{11})^2$ has linear quotient property. Then by Remark \ref{partite} applied to Equation \eqref{eq}, we get that $I(D)^2$ has linear quotient property. \par
  Suppose that $N_D^+(x_{21})\setminus N_D(x_{11})\neq \emptyset$, then there exist a vertex $x_{1t}\in A_1$ such that $x_{1t}\in N_D^+(x_{21})\setminus N_D(x_{11})$  for some $2\leq t\leq n_1$. We have $(x_{21},x_{1t})\in E(D)$ and $x_{11},x_{1t}\in A_1$ then by Remark \ref{partite} we have $x_{1t}\in N_D^+(V^+)$, thus we can write 
  $$I(D)=x_{1t}(x_{21},\ldots,x_{2n_2},\ldots,x_{r1},\dots x_{rn_r})+I(D\setminus x_{1t}).$$
 Therefore, by Lemma \ref{lemma-4.5}, we get that $I(D)^2$ has linear quotient property. 
 \end{proof}

Below we show that if $I(K_{n_1,\dots,n_r})$ is not componentwise linear, then $I(K_{n_1,\dots,n_r})^k$ is not componentwise linear for all $k\geq 2$. 
 
\begin{theorem} \label{thm4}
  Let $D$ be a weighted oriented complete $r$-partite graph, and let \( I(D) \) denote its edge ideal. Then the following statements are equivalent. 
    \begin{enumerate}
        \item[(1)] \( I(D)^k \) is componentwise linear, for some $k\geq 2$; 
        \item[(2)]  \( I(D) \) is componentwise linear; 
        \item[(3)]  $H=H(I(D)_{\ev{2}})$ is co-chordal and $D_i$ are not induced subgraphs of $D$, for $i=1,\ldots,4$. 
    \end{enumerate} 
\end{theorem}
\begin{proof}
$(3)\implies (2)$ follows from Theorem \ref{theorem-3.4}. \par
$(2)\implies (1)$ follows from Theorem \ref{compartite} and Theorem \ref{theorem-4.9}. \par
    $(1)\implies (3)$: Assume $D$ be a $r$-partite weighted oriented graph with $V(D)=A_1\cup \dots\cup A_r$ such that $I(D)^k$ is componentwise linear for some $k\geq 2$. Therefore, by Theorem \ref{D_i}, we have $D_i$ are not induced subgraphs of $D$ for all $i\in\{1,\ldots,4\}$. Since $I(D)^k$ is componentwise linear, then $I(H)^k$ is also componentwise linear because $(I(D)^k)_{\ev{2k}}=I(H)^k$. Then by \cite[Proposition-3.4]{kns25} we get $H^c$ has no induced $4$-cycle. Let $G$ denote the underlying simple graph of $D$. Suppose $H^c$ has an induced cycle $C_m=(x_1,\ldots,x_m)$, for some $m\geq5$. This implies that either the edges $\{x_i,x_{i+1}\}\in E(G)\setminus E(H)$ or the vertices $x_i,x_{i+1}$ lie in same partition in $H$. Then for each $i\in [r]$, we have the induced subgraph $H^c[V(C_m)\cap A_i]$ is an empty graph or a set of isolated vertices or a graph with only one edge. This implies that at least $[\frac{m}{2}]$ edges of $C_m$ are in $E(G)\setminus E(H)$ if $m\geq 6$. And if $m=5$, then is is easy to see that three edges of $C_m$ are in $E(G)\setminus E(H)$. Thus we have at least three edges of $C_m$ are in $E(G)\setminus E(H)$. This gives that at least two vertices (in different partitions), say $x_i,x_j$, incident with these three edges are in $V^+$. Since $\{x_i,x_j\} \notin E(C_m)$ and $C_m$ is an induced subgraph in $H^c$, therefore 
    $\{x_i,x_j\} \notin E(H^c)$. This implies that $\{x_i,x_j\} \in E(H)$. This is a contradiction because $w(x_i) \geq 2$ and $w(x_j) \geq 2$. Therefore $H^c$ has no induced cycle of length $\geq 5$. Hence $H^c$ is chordal. 
\end{proof}

% \begin{question}\label{conjcl}
% Let $D$ be a weighted oriented graph. Are the following equivalent?
% \begin{enumerate}
%     \item $I(D)$ is componentwise linear;
%     \item $I(D)$ has linear quotient property;
%     \item $I(D)$ is a vertex splittable ideal; 
%     \item  $G^c$, $H^c$ is chordal and $D_i$ is not an induced subgraph of $D$ for $i=1,2,3, 4$. 
% \end{enumerate}
% \end{question}

%\noindent 
%After several computational evidences in Macaulay2 \cite{gs}, we raised the following question.
%\begin{question}
% If $D$ is complete $r$-partite, does $I(D)^k$ has linear quotient property, for all $k\geq 2$ if and only if $I(D)$ has linear quotient property? 
%\end{question}

\noindent
{\bf Acknowledgement:} 
Joydip Mondal would like to thank the CSIR-UGC (India) for the financial support through the PhD Fellowship.

% \begin{tikzpicture}[scale=1.5, every node/.style={circle, fill=black, inner sep=1.5pt}]
%   % Define the square vertices
%   \node (A) at (0,0) {};
%   \node (B) at (2,0) {};
%   \node (C) at (2,2) {};
%   \node (D) at (0,2) {};

%   % Define the top triangle vertex
%   \node (E) at (1,3.2) {};

%   % Draw the square
%   \draw (A) -- (B) -- (C) -- (D) -- (A) -- cycle;

%   % Draw the triangle on top
%   \draw (D) -- (E) -- (C);

% \end{tikzpicture}


\begin{thebibliography}{33}
%\bibitem{bdms24} E Basser, R. Diethorn, R. Miranda, M Stinson-Maas,{\em Powers of edge ideals with linear quotients}, arXivprints \url{https://doi.org/10.48550/arXiv.2412.03468}, (2024).

\bibitem{bds23} A. Banerjee, K. Das and S. Selvaraja, {\em Powers of edge ideals of weighted oriented graphs with linear resolutions}, J. Algebra Appl. \textbf{22}, Paper No. 2350148, 12 (2023).

%\bibitem{bed24} Basser, Etan and Diethorn, Rachel and Miranda, Robert and Stinson-Maas, Mario, {\em Powers of Edge Ideals with Linear Quotients},arxivPreprint \url{https://doi.org/10.48550/arXiv.2412.03468}, (2024).

%\bibitem{bh13} S. Bandari, and J. Herzog, {\em Monomial localizations and polymatroidal ideals}, Eur. J. Combin. {\bf 34}, no.~4, 752--763 (2013). 

%\bibitem{chhktt19} G. Caviglia, H. T.  Hà, J. Herzog, M. Kummini, N. Terai and N. Trung, {\em Depth and regularity modulo a principal ideal},  J. Algebraic Combin. \textbf{49}, no. 4, 1-20 (2019).

% \bibitem{ev98} J. Eagon and V. Reiner, {\em Resolutions of Stanley-Reisner rings and Alexander duality},  J. Pure Appl. Algebra. \textbf{130}, no. 3, 265-275 (1998).
 
\bibitem{e95} D. Eisenbud, {\em Commutative algebra with a view toward algebraic geometry}, In: Grad. Texts in Math., vol.150. Springer, New York (1995).

 \bibitem{f90} R. Fr\"{o}berg, {\em On Stanley-Reisner rings}, Topics In Algebra, Part 2 (Warsaw, 1988). \textbf{26, Part 2} pp. 57-70 (1990).

%\bibitem{fh77} S. Foldes and P. Hammer, {\em Split graphs}, Proceedings of the Eighth Southeastern Conference on Combinatorics, Graph Theory and Computing, 311-315 (1977).

\bibitem{fmr24} A. Ficarra, S. Moradi and T. R\"{o}mer, {\em Componentwise linear symbolic powers of edge ideals and Minh’s Conjecture}, arxivPreprint \url{https://doi.org/10.48550/arXiv.2411.11537}, (2024).

 \bibitem{gmsv18} P. Gimenez, J. Mart\'{i}nez-Bernal, A. Simis and R.H. Villarreal, {\em Symbolic powers of monomial ideals and Cohen-Macaulay vertex-weighted digraphs}, in {\it Singularities, algebraic geometry, commutative algebra, and related topics},491--510, Springer, Cham (2018).

\bibitem{gs} D.~R. Grayson and M.~E. Stillman, {\em Macaulay2, a software system for research in algebraic geometry.} Available at \url{http://www.math.uiuc.edu/Macaulay2/}.
 
\bibitem{hh99} J. Herzog and T. Hibi, {\em Componentwise linear ideals}, Nagoya Math. J. \textbf{153} pp. 141-153 (1999).

\bibitem{hhm22} J. Herzog, T. Hibi and S. Moradi, {\em Componentwise linear powers and the x-condition}, Math. Scand. \textbf{128}, no. 3, 401-433 (2022).

%\bibitem{hh11} J. Herzog\ and\ T. Hibi, {\em Monomial ideals}, Graduate Texts in Mathematics, 260, Springer-Verlag London, Ltd., London, (2011).

\bibitem{hhz004} J. Herzog,  T. Hibi, and X. Zheng, {\em Monomial ideals whose powers have a linear resolution}, {\em Math. Scand.} \textbf{95}, 23-32 (2004).

\bibitem{hi05} J. Herzog and S. Iyengar, {\em Koszul modules}, J. Pure Appl. Algebra. \textbf{201}, no. 1-3, 154-188 (2005).

\bibitem{hlmrv19} H. T. H\`a, K. N. Lin, S. Morey, E. Reyes and R. H. Villarreal, {\em Edge ideals of oriented graphs}, Internat. J. Algebra Comput. {\bf 29}, no.~3, 535--559   (2019).

\bibitem{hv22} H.~T. H\`a{} and A. Van~Tuyl, {\em Powers of componentwise linear ideals: the Herzog-Hibi-Ohsugi conjecture and related problems}, Res. Math. Sci. {\bf 9} (2022), no.~2, Paper No. 22, 26 pp.

%\bibitem{ht02} J. Herzog and Y. Takayama, {\em Resolutions by mapping cones},  Homology Homotopy Appl. \textbf{4} pp. 277-294 (2002).

%\bibitem{hhz04} J. Herzog, T. Hibi and X. Zheng, {\em Dirac's theorem on chordal graphs and Alexander duality}, {\em European J. Combin.} \textbf{25}, no. 7, 949-960 (2004).


% \bibitem{kmns13} R. Krithika, R. Mathew, N. Narayanaswamy and, N.  Sadagopan, {\em A Dirac-type characterization of k-chordal graphs }, {\em Discrete Math.} \textbf{313}, no. 24, 2865-2867 (2013).

 \bibitem{kn23} M. Kumar and R. Nanduri, {\em Regularity of powers of edge ideals of Cohen–Macaulay weighted oriented forests},  J. Algebraic Combin. \textbf{58}, no. 3, 867-893 (2023).

 \bibitem{kn25} M. Kumar and R. Nanduri, Regularity of symbolic and ordinary powers of weighted oriented graphs and their upper bounds, Comm. Algebra {\bf 53}, no.~6, 2565--2583 (2025).

\bibitem{kns25} M. Kumar, R. Nanduri, and K. Saha, {\em Componentwise linearity of edge ideals of weighted oriented graphs}, J. Algebraic Combin. \textbf{61(3)}, Paper No. 107881 (2025).

\bibitem{kblo22} S. Kara, J. Biermann, K. N. Lin, and A. O'Keefe, {\em Algebraic invariants of weighted oriented graphs}, J. Algebraic Combin. {\bf 55} , no.~2, 461--491 (2022).

% \bibitem{mka16} S. Moradi and F. Khosh-Ahang, {\em On vertex decomposable simplicial complexes and their Alexander duals}, Math. Scand. \textbf{118}, no. 1, 43-56 (2016).
  
\bibitem{nr15} U. Nagel and T. R\"{o}mer, {\em Criteria for componentwise linearity},  Comm. Algebra. {\bf 43}, no. 3, 935-952 (2015).

\bibitem{prt19} Y. Pitones, E. Reyes\ and\ J. Toledo, {\em Monomial ideals of weighted oriented graphs}, Electron. J. Combin. {\bf 26}, no.~3, Paper No. 3.44, 18 pp  (2019).

\bibitem{reisner76} G. Reisner, {\em Cohen-Macaulay quotients of polynomial rings}, Advances in Math. {\bf 21}, no. 1, 30--49 (1976). 

\bibitem{jz10}  A. Soleyman~Jahan and X. Zheng, {\em Ideals with linear quotients}, J. Combin. Theory Ser. A {\bf 117}, no.~1, 104--110 (2010).
 
%\bibitem{seyed18} S. ~A. Seyed Fakhari, {\em Symbolic powers of cover ideal of very well-covered and bipartite graphs}, Proc. Amer. Math. Soc. {\bf 146}, no. 1, 97--110 (2018).

\bibitem{wg23} H. Wang, G.~J. Zhu and L. Xu, {\em Algebraic properties of powers of edge ideals of vertex-weighted oriented cycles}, Algebra Colloq. {\bf 30} (2023), no.~4, 649--666. 

\bibitem{x21} L. Xu, G. Zhu, H. Wang and J. Zhang, {\em Projective dimension and regularity of powers of edge ideals of vertex-weighted rooted forests}, Bull. Malays. Math. Sci. Soc. {\bf 44}, no.~4, 2215--2233 (2021).  
 

\end{thebibliography}
\end{document}